\documentclass[11pt]{article}
\usepackage[utf8]{inputenc}
\usepackage{amsmath,amssymb,epsfig,bbm}
\usepackage{stmaryrd,mathabx}
\usepackage{comment}
\usepackage{color}
\usepackage[T1]{fontenc}

\usepackage[textsize=small]{todonotes}
\usepackage{enumitem}
\usepackage{varwidth}
\setlist{nolistsep}
\usepackage{hyperref}

\newtheorem{defi}{Definition}[section]
\newtheorem{prop}[defi]{Proposition}
\newtheorem{theo}[defi]{Theorem}
\newtheorem{conj}[defi]{Conjecture}
\newtheorem{lemm}[defi]{Lemma}
\newtheorem{coro}[defi]{Corollary}
\newtheorem{rema}[defi]{Remark}
\newtheorem{exem}[defi]{Example}
\newtheorem{exems}[defi]{Examples}

\newcommand{\bdefi}{\begin{defi}}
\newcommand{\edefi}{\end{defi}}
\newcommand{\bprop}{\begin{prop}}
\newcommand{\eprop}{\end{prop}}
\newcommand{\btheo}{\begin{theo}}
\newcommand{\etheo}{\end{theo}}
\newcommand{\blemm}{\begin{lemm}}
\newcommand{\brema}{\begin{rema}}
\newcommand{\erema}{\end{rema}}
\newcommand{\bexer}{\begin{exem}}
\newcommand{\eexer}{\end{exem}}
\newcommand{\bexems}{\begin{exems}}
\newcommand{\eexems}{\end{exems}}
\newcommand{\bconj}{\begin{conj}}
\newcommand{\econj}{\end{conj}}
\newcommand{\elemm}{\end{lemm}}
\newcommand{\bcoro}{\begin{coro}}
\newcommand{\ecoro}{\end{coro}}
\newcommand{\dem}{\noindent{\bf Proof. }}

\usepackage{mathrsfs}
\renewcommand\mathcal{\mathscr}

\newcommand{\B}{{\cal B}}

\newcommand{\G}{{\cal G}}
\renewcommand{\H}{{\cal H}}
\newcommand{\I}{{\cal I}}

\newcommand{\M}{{\cal M}}
\newcommand{\N}{{\cal N}}
\newcommand{\OOO}{{\cal O}}
\renewcommand{\P}{{\cal P}}

\newcommand{\maths}[1]{{\mathbb #1}}  

\newcommand{\CC}{\maths{C}}

\newcommand{\FF}{\maths{F}}
\newcommand{\HH}{\maths{H}}

\newcommand{\NN}{\maths{N}}

\newcommand{\PP}{\maths{P}}
\newcommand{\QQ}{\maths{Q}}
\newcommand{\RR}{\maths{R}}
\newcommand{\SSS}{\maths{S}}

\newcommand{\ZZ}{\maths{Z}}

\newcommand{\aaa}{{\mathfrak a}}
\newcommand{\bbb}{{\mathfrak b}}

\newcommand{\fff}{{\mathfrak f}}

\newcommand{\mmm}{{\mathfrak m}}

\newcommand{\weakstar}{\overset{*}\rightharpoonup}
\newcommand{\ra}{\rightarrow}
\newcommand{\bs}{\backslash}

\newcommand{\ov}[1]{{\overline #1}} 
\newcommand{\wt}[1]{{\widetilde{#1}}}

\newcommand{\ga}{\gamma}
\newcommand{\Ga}{\Gamma}

\newcommand{\cqfd}{\hfill$\Box$}

\newcommand{\bigO}{\operatorname{O}}
\newcommand{\card}{{\operatorname{Card}}}

\newcommand{\dvol}{\;d\operatorname{vol}}

\newcommand{\haarheis}{\operatorname{Haar}_{\operatorname{Heis}_3}}
\newcommand{\haarhheis}{\operatorname{Haar}_{\operatorname{\HHeis}_7}}

\newcommand{\HHeis}{\HH\!\operatorname{eis}}
\newcommand{\Hom}{\operatorname{Hom}}

\newcommand{\id}{\operatorname{id}}
\renewcommand{\Im}{{\operatorname{Im}}}
\newcommand{\Isom}{\operatorname{Isom}}

\renewcommand{\log}{\operatorname{ln}}

\newcommand{\Par}{\operatorname{Par}}

\newcommand{\Res}{\operatorname{Res}}
\newcommand{\smallo}{\operatorname{o}}

\newcommand{\Vol}{\operatorname{Vol}}
\newcommand{\vol}{\operatorname{vol}}

\newcommand{\hnr}{{\HH}^n_\RR}

\newcommand{\hnh}{{\HH}^n_\HH}

\newcommand{\GL}{\operatorname{GL}}
\newcommand{\Sp}{\operatorname{Sp}}
\newcommand{\PSp}{\operatorname{PSp}}

\newcommand{\PU}{\operatorname{PU}}

\newcommand{\tr}{\operatorname{\tt tr}}
\newcommand{\n}{\operatorname{\tt n}}

\newcommand\normalout{\partial^1_{+}}

\newcommand\normalpm{\partial^1_{\pm}}

\newcounter{fig}

\title{Counting and equidistribution \\
in quaternionic Heisenberg groups}
\author{Jouni Parkkonen \and Fr\'ed\'eric Paulin} 
\begin{document}
\bibliographystyle{/users/jouniparkkonen/Documents/Latex/alphas}
\maketitle
\begin{abstract}
We develop the relationship between quaternionic hyperbolic geometry
and arithmetic counting or equidistribution applications, that arises
from the action of arithmetic groups on quaternionic hyperbolic
spaces, especially in dimension $2$.  We prove a Mertens counting
formula for the rational points over a definite quaternion algebra $A$
over $\QQ$ in the light cone of quaternionic Hermitian forms, as well
as a Neville equidistribution theorem of the set of rational points
over $A$ in quaternionic Heisenberg groups.
\end{abstract}

\section{Introduction}

The two main arithmetic results of this paper are a counting theorem
and an equidistribution theorem of rational points with error
estimates in quaternionic Heisenberg groups, see Theorems
\ref{theo:countintro} and \ref{theo:equidisintro} below.  The proofs
use methods and results from quaternionic hyperbolic geometry,
arithmetic groups and ergodic theory of the geodesic flow in
negatively curved spaces.  We refer for instance to
\cite{Breuillard05, BenQui13a, Kim15} for related results, and
especially to the introductions of \cite{ParPau14AFST,ParPau17MA} for
motivations, going back to the Mertens and Neville counting and
equidistribution results of Farey fractions. The case of the standard
Heisenberg group has been treated in \cite{ParPau17MA}, but new tools
have to be developped in this paper besides dealing with
noncommutativity issues.

Let $\HH$ be Hamilton's quaternion algebra over $\RR$, with $x\mapsto
\overline{x}$ its conjugation, $\n: x\mapsto x\overline{x}$ its
reduced norm, $\tr: x\mapsto x+\overline{x}$ its reduced trace. Let
$A$ be a definite ($A\otimes_\QQ\RR=\HH$) quaternion algebra over
$\QQ$, with discriminant $D_A$ and class number $h_A$. Let $m_A=72$ if
$D_A$ is even, and $m_A=1$ otherwise. Let $\OOO$ be a maximal order in
$A$. We denote by $_\OOO\langle a,\alpha,c\rangle$ the left ideal of
$\OOO$ generated by $a,\alpha,c\in \OOO$.  See \cite{Vigneras80} and
Section \ref{sect:quaternion} for definitions. The $2$-step nilpotent
group
\begin{equation}\label{eq:heisintro}
\N(\OOO)=\{(w_0,w)\in \OOO\times \OOO\;:\; \tr(w_0)=\n(w)\}
\end{equation} 
with law 
\begin{equation}\label{eq:heislaw}
(w_0,w)(w'_0,w')= (w_0+w'_0+\,\overline{w}\,w',w+w')
\end{equation} 
acts on $\OOO\times\OOO\times\OOO$ by the shears
\begin{equation}\label{eq:heisaction}
(w_0,w)(a,\alpha,c)=(a+\overline{w}\,\alpha+w_0\,c,\alpha+w\,c,c)\;.
\end{equation}

\btheo\label{theo:countintro} There exists $\kappa>0$
such that as $s\ra+\infty$,
\begin{align*} 
\card\;\;\N(\OOO)\bs
\big\{(a,\alpha,c)&\in \OOO\times \OOO  \times\OOO\;:\;\, 
_\OOO\langle a,\alpha,c\rangle =\OOO, \;\tr(\overline{a}\,c)=
\n(\alpha), \;\n(c)\leq s\big\} \\ & =\;  
\frac{204\,120\;D_A^{\;4}}
{\pi^8\;m_A\;|\OOO^\times|\;\prod_{p|D_A}(p-1)(p^2+1)(p^3-1)} 
\;s^5\,(1+\bigO(s^{-\kappa}))\,. 
\end{align*}
\etheo

The quaternionic Heisenberg group
$$
\HHeis_{7}= \{(w_0,w)\in\HH\times\HH \;:\; \tr\;w_0=\n(w)\}\,,
$$ 
with the group law defined by Equation \eqref{eq:heislaw} is the Lie
group of $\RR$-points of a $\QQ$-group whose group of $\QQ$-points is
$\HHeis_{7}\cap (A\times A)$, and in which $\N(\OOO)=\HHeis_{7}\cap
(\OOO\times \OOO)$ is a (uniform) lattice. We endow it with its Haar
measure $\haarhheis$ normalised in such a way that the total mass of
the induced measure on $\N(\OOO)\bs\HHeis_{7}$ is $\frac{D_A^2}{4}$.
We will explain later on this normalisation.  Theorem
\ref{theo:countintro} is a counting result of rational points
$(ac^{-1},\alpha c^{-1})$ (analogous to Farey fractions) in
$\HHeis_{7}$, and the following result is a related equidistribution
theorem. In this paper, we denote by $\Delta_x$ the unit Dirac mass at
a point $x$.

\btheo\label{theo:equidisintro} As $s\ra+\infty$, we have
\begin{align*}
\frac{\pi^8\;m_A\;|\OOO^\times|^2\;\prod_{p|D_A}(p-1)(p^2+1)(p^3-1)}
{816\,480\;D_A^2}\;s^{-5}&\;\times \\
\;\sum_{\substack{(a,\,\alpha,\,c)\in \OOO\times\OOO\times\OOO,
\;0<\n(c)\leq s \\ \tr(a\,\overline{c})=\n(\alpha),
\;_\OOO\langle a,\,\alpha,\,c\rangle=\OOO}}\;
&\Delta_{(ac^{-1},\,\alpha c^{-1})}\;\weakstar\;\haarhheis\;.
\end{align*}
\etheo

We refer to Theorems \ref{theo:countHeis} and \ref{theo:equidisHeis}
in Section \ref{sect:mertens} for more general results with added
congruence properties, and to Remark \ref{rem:highdim} for counting
and equidistribution results in higher dimensional quaternionic
Heisenberg groups.
%
%

\medskip The proof of the above arithmetic results strongly rely on
quaternionic hyperbolic geometry that we recall and develop in
Sections \ref{sec:quathypspace} and \ref{sec:quathypgeom} (see also
for instance \cite{Allcock00,KimPar03,CaoPar11,Kim15,Philippe16,
  CaoPar18,EmeKim18}). Let $q$ be the quaternionic Hermitian form on
$\HH^{3}$ defined by
%
%
$$
q(z_0,z_1, z_2)=-\tr(\,\overline{z_0}\,z_2) + \n(z_1)\;,
$$ 
and $\PU_q$ its projective unitary group, which is the isometry
group of the quaternionic hyperbolic plane $\HH^2_\HH$, realized as
the negative cone of $q$ in the right projective plane $\PP^2_{\rm
  r}(\HH)$, and normalised in order to have maximal sectional
curvature $-1$. The proofs of Theorems \ref{theo:countintro} and
\ref{theo:equidisintro} use the following two results of independent
interests.

The subgroup $\PU_q(\OOO)=P(\GL_{3}(\OOO)\cap U_q)$ of $\PU_q$ is an
arithmetic lattice, and hence by a standard result of
Borel-Harish-Chandra, the orbifold $\PU_q(\OOO)\bs \HH^2_\HH$ has
finitely many ends (also called {\it cusps}). By adapting Zink's
method \cite{Zink79}, we compute their number in Section
\ref{sec:cuspnumb}.

\btheo\label{intro:cuspnumb} The number of ends of the quaternionic
hyperbolic orbifold $\PU_q(\OOO)\bs \HH^2_\HH$ is equal to the class
number $h_A$ of $A$.  
\etheo

Using Prasad's formula \cite{Prasad89} and adapting Emery's Appendix
in \cite{ParPau13ANT}, we compute in Section \ref{sec:computvol} the
volume of the quaternionic hyperbolic manifolds $\PU_q(\OOO) \bs
\HH^2_\HH$.  This computation (equivalent by Hirzebruch's
proportionality theorem to the computation of its  Euler-Poincaré
characteristic), is close to, and uses argument from, the paper
\cite{EmeKim18} of Emery-Kim.

\btheo\label{intro:computvol} 
%
%
The volume of the quaternionic hyperbolic orbifold $\PU_q(\OOO)\bs
\HH^2_\HH$ is equal to
$$
\Vol(\PU_q(\OOO_q)\bs \HH^2_\HH)=
\frac{\pi^4\;m_A}{42525\;2^{13}}\;\prod_{p|D_A}(p-1)(p^2+1)(p^3-1)\;.
$$  
\etheo

The analytical tools for the proof of Theorems \ref{theo:countintro}
and \ref{theo:equidisintro} are developped in Section
\ref{sect:measurecomput}, where we give various computations of the
measures that appear in the application of the ergodic tools of
\cite{ParPau17ETDS}, see also \cite[Chap.~12]{BroParPau19} where these
results are announced.

The Cygan distance on the quaternionic Heisenberg group, the Poisson
kernel, the Patterson measures introduced and computed in Section
\ref{sect:measurecomput}, and related quantities, should be useful in
potential theory on the quaternionic Heisenberg group and for the
study of the hypoelliptic Laplacian in sub-Riemannian geometry, see
for instance \cite{FolSte74,Krantz09} for the (complex) Heisenberg
group. 

In the subsequent paper \cite{ParPau20b}, we will give geometrical
applications of this paper to counting and equidistribution of
quaternionic chains in the boundary at infinity of the quaternionic
hyperbolic plane. A quaternionic chain is the boundary at infinity of a
quaternionic geodesic line (as defined in Section
\ref{sec:quathypspace}). We will also prove a Cartan-type theorem of
rigidity for the bijections of the $\partial_\infty\hnh$ preserving the
set of quaternionic chains.

\medskip\noindent{\small {\it Acknowledgements: } We thank Gaetan
Chenevier for his help with Section \ref{sec:cuspnumb}. The second
author thanks Hee Oh and the hospitality of Yale University for a
one month stay, decisive for the writing of this paper.}

\section{A reminder on quaternion algebras}
\label{sect:quaternion}

In this section, we recall the basic definitions and a few facts on
quaternion algebras ($4$-dimensional central simple algebras),
quaternionic linear algebra, and the ideal theory in maximal orders in
quaternion algebras. Our main reference for this material is
\cite{Vigneras80}.

Given a skew field $D$ and $n\in\NN$, we denote by $\PP^n_{\rm r}(D)$
the right\footnote{For all $x_0,x_1,\dots, x_N\in D$, we have
  $(x_1,\dots,x_N)x_0= (x_1x_0,\dots,x_Nx_0)$.} projective space of
$D$ of dimension $n$, that is, the space of lines of the right vector
space $D^{n+1}$ over $D$.

Let $\HH$ be Hamilton's quaternion algebra over $\RR$, with $x\mapsto
\overline{x}$ its conjugation, $\n: x\mapsto x\overline{x}$ its
reduced norm, $\tr: x\mapsto x+\overline{x}$ its reduced trace. Note
that $\n(xy)=\n(x)\n(y)$, $\tr(\overline{x})=\tr(x)$ and
$\tr(xy)=\tr(yx)$ for all $x,y\in\HH$. Let
$$
\Im\;\HH=\{x\in\HH\;:\;\tr x=0\}
$$ 
be the $\RR$-subspace of purely imaginary quaternions of $\HH$, so
that every $x\in \Im\;\HH$ satisfies $\overline{x}=-x$. For every
$x\in\HH$, let $\Im\;x=x-\frac{1}{2}\,\tr(x)=\frac 12(x-\overline x)$,
which is a purely imaginary quaternion.

For every $N\in \NN-\{0\}$, we consider the right vector space $\HH^N$
over $\HH$, on which the group $\GL_N(\HH)$ acts linearly on the left.
For all $w=(w_1,\dots,w_{N})$ and $w'=(w'_1,\dots,w'_{N})$ in
$\HH^{N}$, we denote by $\overline{w}\cdot w'= \sum_{p=1}^{N}
\overline{w_p}\,w'_p$ their standard quaternionic Hermitian
product,\footnote{We have $\overline{w\lambda}\cdot (w'\mu)=
  \overline{\lambda}\, (\overline{w}\cdot w')\mu$ for all
  $w,w'\in\HH^{N}$ and $\lambda,\mu\in\HH$.} and we define $\n(w)=
\overline w\cdot w= \sum_{p=1}^{N} \n(w_p)$. We endow $\HH^{N}$ with
the standard Euclidean structure $(w,w')\mapsto \frac12
\tr(\,\overline{w}\cdot w')$. In particular, $\HH$ and $\Im\;\HH$ are
endowed with the Euclidean structure making their standard basis
$(1,i,j,k)$ and $(i,j,k)$ orthonormal.

\medskip
On the right vector space $\HH\times\HH^{n-1}\times\HH$ over $\HH$
with coordinates $(z_0,z,z_n)$, let $q$ be the nondegenerate
quaternionic Hermitian form\footnote{It does satisfy
  $q(x\lambda)=\n(\lambda)\,q(x)$ for all $\lambda \in \HH$ and
  $x\in\HH^{n+1}$, since $\tr(uv)=\tr(vu)$ for all $u,v\in \HH$.}
\begin{equation}\label{eq:q}
q(z_0,z,z_n)=-\tr(\,\overline{z_0}\,z_n) + \n(z)
\end{equation}
of Witt signature $(1,n)$, and let $\Phi:\HH^{n+1}\times \HH^{n+1}\ra
\HH$, defined by
\begin{equation}\label{eq:Phi}
\Phi:((z_0,z,z_n),(z'_0,z',z'_n))\mapsto
-\overline{z_0}\,z'_n-\overline{z_n}\,z'_0+\,\overline{z}\cdot z'\,,
\end{equation}
be the associated quaternionic sesquilinear form.  An element $x\in
\HH^{n+1}$ is {\it isotropic} if $q(x)=0$.
%
%

\bigskip
Throughout this paper, $A$ is a quaternion algebra over $\QQ$, which
is {\em definite}, that is, $A\otimes_\QQ\RR$ is isomorphic with
$\HH$.  We fix an identification of $A\otimes_\QQ\RR$ and $\HH$ and,
accordingly, consider $A$ as a $\QQ$-subalgebra of $\HH$.

The {\em reduced discriminant} $D_A$ of $A$ is the product of the
primes $p\in\NN$ such that $A\otimes_\QQ\QQ_p$ is a division algebra.
Two definite quaternion algebras over $\QQ$ are isomorphic if and only
if they have the same reduced discriminant, which can be any product
of an odd number of primes (see \cite[page 74]{Vigneras80}).

A {\em $\ZZ$-lattice} $I$ in $A$ is a finitely generated
$\ZZ$-submodule of $A$ generating $A$ as a $\QQ$-vector space. An {\it
  order} in $A$ is a unitary subring $\OOO$ of $A$ which is a
$\ZZ$-lattice. It is contained in a {\it maximal order} (for the
inclusion).  The {\em left order} of a $\ZZ$-lattice $I$ is
$$\OOO_\ell(I)=\{x\in A\;:\; xI \subset I\}\,.$$ 

From now on, let $\OOO$ be a maximal order in $A$. It is well known
that the trace map $\tr:\OOO\ra \ZZ$ is surjective.\footnote{See, for
  instance, the proof of Prop.~16 in \cite{ChePau19}.} A {\em left
  fractional ideal} of $\OOO$ is a $\ZZ$-lattice of $A$ whose left
order is $\OOO$. A {\em left (integral) ideal} of $\OOO$ is a left
fractional ideal of $\OOO$ contained in $\OOO$.  For any subset $B$ of
$A$, we denote by $_\OOO\langle B\rangle$ the left fractional ideal of
$\OOO$ generated by the elements of $B$.  Right fractional ideals are
defined analogously. The {\em inverse} of a left fractional ideal
$\mmm$ of $\OOO$ is the right fractional ideal
$$
\mmm^{-1}=\{x\in A\;:\;\mmm \,x \,\mmm\subset \mmm\}=
\{x\in A\;:\;\mmm \,x \subset \OOO\}\;.
$$
For all $u,v\in\OOO-\{0\}$, we have
\begin{equation}\label{eq:invsomeginter}
(\OOO u+ \OOO v)^{-1}=u^{-1} \OOO\cap v^{-1}\OOO\;.
\end{equation}

If $M$ is a right $\OOO$-module, then endowed with the pointwise
multiplication by $\OOO$ on the left, the $\ZZ$-module
$\Hom_\OOO(M,\OOO)$ (of morphisms of right $\OOO$-modules from $M$ to
$\OOO$) is a left $\OOO$-module. We denote by $\widecheck{M}$ the left
$\OOO$-module equal to the $\ZZ$-module $M$ endowed with the left
multiplication by $\OOO$ defined by $(\lambda,v)\mapsto
v\,\overline{\lambda}$. If $\mmm$ is a right fractional ideal of
$\OOO$, then the map from $\mmm^{-1}$ to $\Hom_\OOO(\mmm,\OOO)$
defined by $x\mapsto \{y\mapsto xy\}$ is an isomorphism of left
$\OOO$-modules, see for instance \cite[page 192]{Reiner75}.

Two left fractional ideals $\mmm$ and $\mmm'$ of $\OOO$ are isomorphic
as left $\OOO$-modules if and only if $\mmm'=\mmm c$ for some $c\in
A^\times$. A {\it (left) ideal class} of $\OOO$ is an equivalence
class of left fractional ideals of $\OOO$ for this equivalence
relation. We will denote by $_\OOO\!\I$ the set of ideal classes of
$\OOO$.
The {\it class number} $h_A$ of $A$ is the number of ideal classes of
$\OOO$. It is finite and independent of the maximal order $\OOO$.

We denote by $\OOO^\times$ the group of invertible elements (or
equivalently of norm $1$ elements) of $\OOO$. Its order is $2$, $4$ or
$6$ except that $|\OOO^\times|=24$ when $D_A=2$ and $|\OOO^\times|=12$
when $D_A=3$ (see \cite[page 103]{Eichler38} for a formula when
$h_A=1$).

By for instance \cite[Lem.~5.5]{KraOse90}, the covolume of the
$\ZZ$-lattice $\OOO$ in the Euclidean vector space $\HH$ is
\begin{equation}\label{eq:covolmaxorder}
\Vol(\OOO\bs\HH)=\frac{D_A}{4}\;.
\end{equation}

\section{Quaternionic hyperbolic space}
\label{sec:quathypspace}

In this section, we recall some background on the quaternionic
hyperbolic spaces, as mostly contained in \cite{KimPar03}, see also
\cite{Philippe16}. Note, however, that our conventions differ from
those of these references in the sesquilinearity properties of
Hermitian products, in the choice of the Hermitian form of Witt
signature $(1,n)$, and in the normalisation of the curvature.

We fix $n\in\NN-\{0,1\}$.  The {\it Siegel domain} model of the
quaternionic hyperbolic $n$-space $\hnh$ is
$$ 
\big\{(w_0,w)\in\HH\times\HH^{n-1}\;:\; 
\tr\, w_0 -\n(w)>0\big\}\,,
$$ 
endowed with the Riemannian metric
\begin{equation}\label{eq:Siegdommet}
ds^2_{\,\hnh}=\frac{1}{(\tr\, w_0 -\n(w))^2}
\big(\,\n(dw_0-\overline{dw}\cdot w)+
(\tr\, w_0-\n(w))\;\n(dw)\,\big)\,.
\end{equation}
Note that this metric is normalised so that its sectional curvatures
are in $[-4,-1]$, instead of in $[-1,-\frac{1}{4}]$ as in
\cite{KimPar03} and \cite{Philippe16}. This will facilitate in Section
\ref{sect:mertens} the references to works using that
normalisation. Its boundary at infinity is
$$
\partial_\infty\hnh=\big\{(w_0,w)\in
\HH\times \HH^{n-1} \;:\; \tr\, w_0 -\n(w)=0\big\}\cup\{\infty\}\,.
$$

A {\it quaternionic geodesic line} in $\hnh$ is the image by an
isometry of $\hnh$ of the intersection of $\hnh$ with the quaternionic
line $\HH\times\{0\}$. With our normalisation of the metric, a
quaternionic geodesic line is a totally geodesic submanifold of real
dimension $4$ and constant sectional curvature $-4$.

The Siegel domain $\hnh$ embeds in the right quaternionic projective
$n$-space $\PP^n_{\rm r}(\HH)$ by the map (using homogeneous
coordinates)
$$
(w_0,w)\mapsto [w_0:w:1]\;.
$$ 
By this map, we identify $\hnh$ with its image, which when endowed
with the isometric Riemannian metric, is called the {\it projective
  model} of $\hnh$. Note that this image is the {\it negative cone} of
the quaternionic Hermitian form $q$ defined in Equation \eqref{eq:q},
that is $\big\{[z_0:z:z_n]\in \PP^n_{\rm r}(\HH) \;:\;q(z_0,z,z_n)
<0\big\}$.  This embedding extends continuously to the boundary at
infinity, by mapping the point $(w_0,w)\in \partial_\infty\hnh
-\{\infty\}$ to $[w_0:w:1]$ and $\infty$ to $[1:0:0]$, so that the
image of $\partial_\infty\hnh$ is the {\it isotropic cone} of $q$,
that is $\big\{[z_0:z:z_n]\in\PP^n_{\rm r}(\HH)\;:\; q(z_0,z,z_n)=0
\big\}$.

The distance between two points in the Siegel domain of $\hnh$ has an
explicit expression using the projective model: If $(w_0,w),
(w_0',w')\in\hnh$, with $\Phi$ defined in Equation \eqref{eq:Phi},
then
\begin{equation}\label{ed:distexplicit}
\cosh^2 d((w_0,w),(w_0',w'))=
\frac{\Phi((w_0,w,1),(w_0',w',1))\; \Phi((w_0',w',1),(w_0,w,1))}
{q(w_0,w,1)\;q(w_0',w',1)}\,,
\end{equation}
see for example
\cite{Mostow73} with the same normalisation of the metric as ours,
\cite[page 292]{KimPar03} and \cite[Sect.~1.2]{Philippe16} with a
discussion of the different normalizations of the curvature.

For every $N\in\NN$, let $I_N$ be the identity $N\times N$ matrix. Let
$$
J=\begin{pmatrix} 0 & 0 & -1\\0 & I_{n-1} & 0\\-1 & 0 & 0
\end{pmatrix}\;,
$$ 
which differs only up to signs with the matrix $J$ in \cite{KimPar03}.
Given a quaternionic matrix $X=(x_{p,p'})_{1\leq p\leq r,\, 1\leq p'\leq
  s} \in \M_{r,s}(\HH)$, we denote by $X^*= (x^*_{p,p'}=
\overline{x_{p',p}}\,)_{1\leq p\leq s,\,1\leq p'\leq r}\in
\M_{s,r}(\HH)$ its {\it conjugate-transpose} matrix.  Let
$$
\operatorname{U}_q=\{g\in \GL_{n+1}(\HH)\;:\; q\circ g=q\}
=\{g\in \GL_{n+1}(\HH)\;:\; g^*J\,g=J\}
$$ 
be the {\it unitary group} of $q$.  Its linear action on $\HH^{n+1}$
induces a projective action on $\PP^n_{\rm r}(\HH)$ with kernel its
center, which is reduced to $\{\pm I_{n+1}\}$.  The {\it projective
  unitary group}
$$
\PU_q=\operatorname{U}_q/\{\pm I_{n+1}\}
$$ 
of $q$ acts faithfully on $\PP^n_{\rm r}(\HH)$, preserving $\hnh$, and
its restriction to $\hnh$ is the full isometry group of $\hnh$.  The
connected (almost-)simple real Lie groups $\operatorname{U}_q$ and
$\PU_q$ are also denoted by $\Sp(1,n)$ and $\PSp(1,n)$, when the
dependence on the choice of $q$ is not important. 

We identify any
element of $\HH^{n-1}$ with its column matrix. If
$$
X=\begin{pmatrix} a & \gamma^* & b\\ \alpha & A & \beta 
\\ c & \delta^* & d\end{pmatrix}\in \GL_{n+1}(\HH)
$$ 
is a matrix with $a,b,c,d\in\HH$, $\alpha,\beta,\gamma,\delta
\in\HH^{n-1}$ and $A\in\M_{n-1,n-1}(\HH)$, then
$$
JX^*J=\begin{pmatrix} \overline{d} & -\beta^* &
\overline{b}\\ -\delta & A^* & -\gamma \\ \overline{c} & -\alpha^*
& \overline{a}\end{pmatrix}\,.
$$ 
The matrix $X$ belongs to $\operatorname{U}_{q}$ if and only if $X$ is
invertible with inverse $JX^*J$.  In particular,  $X$ belongs to
$\operatorname{U}_{q}$ if and only if 
\begin{equation}\label{eq:equationsUq}
\begin{cases}
\hfill c\,\overline{d}-\delta^*\delta+d\,\overline{c}   &= 0 \\ 
\hfill  a\,\overline{b}-\ga^*\ga+b\,\overline{a}   &= 0 \\ 
   -\alpha\beta^*+A A^* - \beta\alpha^*  \!\!\!&=I_{n-1} \\
\hfill  c\,\overline{b}-\delta^*\ga+d\,\overline{a} &= 1\\
\hfill  \alpha \,\overline{d}-A\delta+\beta\,\overline{c} &=0\\
\hfill  \alpha\,\overline{b}-A\ga+\beta\,\overline{a}  &=0 \;.
\end{cases}
\end{equation}
These equations are the same ones as in the complex hyperbolic case in
\cite[\S6.1]{ParPau10GT}, up to being careful with the orders of the
products; see also \cite{CaoPar11,CaoPar18} with different sign
conventions.

By for instance \cite{KimPar03} or the set of equations
\eqref{eq:equationsUq}, an element $g\in \operatorname{U}_q$ fixes
$\infty$ if and only if its $(1,3)$ entry vanishes, or, equivalently,
if $g$ is block upper triangular (this is the reason, besides
rationality problems, that we chose the quaternionic Hermitian form
$q$ rather than a diagonal one).  We denote by
$$
\Sp(1)=\{x\in\HH\;:\;\n(x)=1\}
$$ 
the subgroup of units of norm one of $\HH^\times$, and
$$
\Sp(n-1)= \{g\in \GL_{n-1}(\HH)\;:\; g^*g=I_{n-1}\}
$$ 
the compact symplectic group in dimension $n-1$.  An easy computation
shows that the block upper triangular subgroup of $\operatorname{U}_q$
is
$$
\operatorname{B}_q=\Bigg\{\begin{pmatrix} 
\mu r & \zeta^* & \frac{1}{2r}(\n(\zeta)+u)\mu \\
0 & U & \frac{1}{r}\;U\,\zeta\;\mu\\
0 & 0 & \frac{\mu}{r}
\end{pmatrix}\;:\; 
\begin{array}{c}\zeta\in \HH^{n-1},\; u\in\Im\;\HH,\\
U\in \Sp(n-1), \mu\in \Sp(1), r>0\end{array}\Bigg\}\;.
$$ 
Its image $\operatorname{PB}_q=\operatorname{B}_q/\{\pm I_{n+1}\}$
in $\PU_q$ is equal to the stabiliser of $\infty$ in $\PU_q$.

\section{The number of cusps of $\PU_q(\OOO)$}
\label{sec:cuspnumb}
%
%

Let $A$ be a definite quaternion algebra over $\QQ$, and let $\OOO$ be
a maximal order in $A$. Let $q$ be the quaternionic Hermitian form
defined in Equation \eqref{eq:q}.  Let $\operatorname{U}_q (\OOO) =
\operatorname{U}_q\cap\GL_{n+1}(\OOO)$, which is (see below) an
arithmetic lattice in $\operatorname{U}_q$, and let $\PU_q(\OOO)$ be
its image in $\PU_q$. The aim of this section is to describe precisely
the structure of the set of ends of the finite volume quaternionic
hyperbolic orbifold $\PU_q(\OOO)\bs\HH_\HH^n$ when $n=2$.

In this section and the following one, we will need to make explicit
the arithmetic structure of $\operatorname{U}_q(\OOO)$.  Since $J$ has
rational coefficients, we consider the linear algebraic group
$\underline{G}$ defined over $\QQ$, such that $\underline{G}(\QQ)
=\{g\in \GL_{n+1}(A)\;:\; g^*J\,g=J\}$, and $\underline{G}(K)=\{g\in
\GL_{n+1}(A\otimes_\QQ K)\;:\; g^*J\,g=J\}$ for every commutative
field $K$ with characteristic $0$. In particular, $\underline{G}(\RR)=
\operatorname{U}_q= \Sp(1,n)$ and $\underline{G}(\CC) \simeq
\Sp_{n+1}(\CC)$.\footnote{This group is sometimes denoted by
  $\Sp_{2(n+1)}(\CC)$, for instance in \cite{PlaRap94}.} By \cite[\S
  2.3.3]{PlaRap94}, the algebraic group $\underline{G}$ is absolutely
connected, (quasi-)simple, simply connected of type $C_{n+1}$, and
denoted by $U_{n+1}(A,q)$ in loc.~cit.

Considering $\OOO^{n+1}$ as a $\ZZ$-lattice of $\HH^{n+1}$, we endow
$\underline{G}$ with the natural $\ZZ$-form such that
$\underline{G}(\ZZ) =\operatorname{U}_q (\OOO)$ and
$\underline{G}(\ZZ_p) =\{g\in \GL_{n+1}(\OOO_p)\;:\; g^*Jg=J\}$ for every
prime $p$, where $\OOO_p=\OOO\otimes_\ZZ\ZZ_p$.

Let us recall a few facts that follow from the work of Borel and
Harish-Chandra (see for instance \cite[Th.~1.10]{Borel66b}). The
discrete group $\underline{G}(\ZZ)$ is a lattice in $\underline{G}
(\RR)$. If $\underline{P}$ is a minimal parabolic subgroup of
$\underline{G}$ defined over $\QQ$ (for instance the stabiliser of
$\infty$), then the set $\operatorname{Par}_{q,\,\OOO}$ of parabolic
fixed points of $\underline{G}(\ZZ)$ in $\underline{G}(\RR)/
\underline{P}(\RR) = \partial_\infty\HH^n_\HH$ is exactly
$$
\operatorname{Par}_{q,\,\OOO} = \underline{G}(\QQ)
\underline{P}(\RR) = \partial_\infty\HH^n_\HH\cap \PP_{\rm r}^n(A)\;,
$$
This is the set of isotropic rational projective points in $\PP_{\rm
  r}^n(\HH)$, on which $\underline{G}(\ZZ)$ acts with finitely many
orbits. In particular, the {\it set of cusps} $\PU_q(\OOO)\bs
\Par_{q,\OOO}$ is in bijection with $\underline{G}(\ZZ)\bs
\underline{G}(\QQ)/\underline{P}(\QQ)$.

For every right $\OOO$-submodule $M$ of $\OOO^{n+1}$, with
$\Phi_A:A^{n+1}\times A^{n+1}\ra A$ the restriction over $A$ of the
form $\Phi$ defined in Equation \eqref{eq:Phi}, we denote by
$$
M^\perp =\{y\in\OOO^{n+1}\;:\;\forall\;x\in M,\;\Phi_A(x,y)=0\}
$$  
the right $\OOO$-submodule of $\OOO^{n+1}$ orthogonal to $M$.  Note
that $\Phi_A(\OOO^{n+1}\times\OOO^{n+1})= \OOO$.  The Hermitian
$\OOO$-module $(\OOO^{n+1},\Phi_A)$ is {\it unimodular}, that is, the
map
\begin{align}\label{eq:theta} 
\Theta:\;\widecheck{\OOO^{n+1}}\; &\ra \Hom_\OOO(\OOO^{n+1},\OOO)\nonumber\\ 
z \;& \mapsto \{z'\mapsto \Phi_A(z,z')\}
\end{align}
is an isomorphism of left $\OOO$-modules. It is indeed clearly an
injective morphism of left $\OOO$-modules. Its surjectivity comes from
the fact that the coordinate forms $z'\mapsto z'_n$, $z'\mapsto z'_0$,
$z'\mapsto z'_i$ for $1\leq i\leq n-1$ are up to signs the images by
$\Theta$ of the canonical basis elements $e_0$, $e_n$ and $e_i$ for
$1\leq i\leq n-1$ respectively.

For every $x=(x_0,x_1,\dots, x_n)\in A^{n+1}$, let $_\OOO\langle
x\rangle=\OOO x_0+\OOO x_1+\dots+ \OOO x_n$ be the left fractional
ideal of $\OOO$ generated by $x_0,x_1,\dots, x_n$. The proof of the following
result adapts arguments of \cite{Zink79}, where $\OOO$ is replaced by
the ring of integers of an imaginary quadratic field.

\bprop\label{prop:isotropictransitive} Assume that $n=2$. For all
isotropic elements $x,x'\in A^{n+1}$, there exists an element in
$\PU_q(\OOO)$ sending the image of $x$ in $\PP^n_{\rm r}(A)$ to the
one of $x'$ if and only if the left fractional ideals $_\OOO\langle
x\rangle$ and $_\OOO\langle x'\rangle$ have the same class.  
\eprop

We do not know whether the result remains valid when $n\geq 3$.

\medskip
\dem The direct implication is immediate. Conversely, let
$x=(x_0,\dots, x_n)$ and $x'$ be isotropic elements of
$\OOO^{n+1}$ such that $_\OOO\langle x\rangle$ and $_\OOO\langle
x'\rangle$ are in the same left ideal class. This means that there is
some $c\in A^\times$ such that $_\OOO\langle x\rangle =\;_\OOO\langle
x'\rangle c=\;_\OOO\langle x'c\rangle$.  In particular, this implies
that $x'c\in \OOO^{n+1}$. As we are interested in the images of $x$
and $x'$ in $\PP^n_{\rm r}(A)$, it is therefore sufficient to prove
that there exists an element of $\operatorname{U}_q(\OOO)$ sending $x$
to $x'$ if $_\OOO\langle x\rangle =\;_\OOO\langle x'\rangle$. 

Let $\aaa=\{a\in A\;:\; xa\in\OOO^{n+1}\}$ and $\aaa'=\{a\in A\;:\;
x'a\in\OOO^{n+1}\}$, which are right fractional ideals of $\OOO$,
containing $1$ since $x,x'\in\OOO^{n+1}$.  By Equation
\eqref{eq:invsomeginter}, we have (omitting $x_i^{\;-1}\OOO$ and $\OOO
x_i$ if $x_i=0$)
$$
\aaa=x_0^{\;-1}\OOO\cap\dots\cap x_n^{\;-1}\OOO= 
(\OOO x_0+\dots+ \OOO x_n)^{-1}=(\,_\OOO\langle x\rangle\,)^{-1}
=(\,_\OOO\langle x'\rangle\,)^{-1}=\aaa'\;.
$$

Composing the map $\Theta$ defined in Equation \eqref{eq:theta} with
the restriction map to $x\aaa$, we have a surjective morphism of left
$\OOO$-modules from $\widecheck{\OOO^{n+1}}$ to
$\Hom_\OOO(x\aaa,\OOO)$. Its kernel is the orthogonal subspace
$(x\aaa)^\perp =\{y\in\OOO^{n+1}\;:\; \Phi_A(x,y)=0\}$, which contains
$x\aaa$ since $x$ is isotropic.

Let $y\in A$ be such that $\Phi_A(x,y) \neq 0$, which exists since
$\Phi_A$ is nondegenerate. Up to replacing $y$ by $y\,\Phi_A(x,y)^{-1}$,
we may assume that $\Phi_A(x,y)=1$. Let $\mmm$ be the right fractional
ideal of $\OOO$ such that $\OOO^{n+1}=(x\aaa)^\perp\oplus
y\mmm$. Composing by the explicit isomorphisms of right $\OOO$-modules
$$
\mmm\simeq y\mmm\simeq \OOO^{n+1}/(x\aaa)^\perp\simeq
\widecheck{\Hom}_\OOO(x\aaa,\OOO)\simeq
\widecheck{\Hom}_\OOO(\aaa,\OOO)\simeq\widecheck{\aaa}^{\;-1}\;,
$$ 
we have $\mmm=\widecheck{\aaa}^{\;-1}$. Hence there exists a right
$\OOO$-submodule $M$ of $\OOO^{n+1}$ such that
$$
\OOO^{n+1}=x\aaa \oplus y \,\widecheck{\aaa}^{\;-1}\oplus M\;.
$$ 
Note that the map $\tr:A\ra \QQ$ is onto. Since $\Phi_A(y+ x\lambda,y+
x\lambda)= \Phi_A(y,y)+\tr\lambda$, up to replacing $y$ by $y+ x\lambda$
for some $\lambda\in A$ such that $\tr\lambda= -\Phi_A(y,y)$, which is
possible since $\Phi_A(y,y)\in A\cap\RR=\QQ$, we may assume that
$q(x)=q(y)=0$ and $\Phi_A(x,y)=1$.

Since $x\aaa \oplus y \,\widecheck{\aaa}^{\;-1}$ is unimodular, we may
take $M=(x\aaa\oplus y\, \widecheck{\aaa}^{\;-1})^\perp$. Since $n=2$,
we may write $M= z\bbb$ for some $z\in A^{n+1}$ such that
$\Phi_A(x,z)=\Phi_A(y,z)=0$ and some right fractional ideal $\bbb$ of
$\OOO$. Since $(\OOO^{n+1},\Phi_A)$ is unimodular and $z\bbb$ is
orthogonal to $x\aaa \oplus y \,\widecheck{\aaa}^{\;-1}$, we have
$\Phi_A(z\bbb,z\bbb)=\overline{\bbb}\,\bbb \,q(z)$ contains $1$ and is
contained in $\OOO$, hence is equal to $\OOO$. Therefore
$q(z)=\frac{1}{\n(\bbb)}$.

Similarly, we have $\OOO^{n+1} = x'\aaa \oplus y'
\,\widecheck{\aaa}^{\;-1} \oplus z'\bbb'$ with
$$
q(x')=q(y')=\Phi_A(x',z')=\Phi_A(y',z')=0\;\;\;{\rm and}\;\;\;
q(z')=\frac{1}{\n(\bbb')}\;.
$$ 
Since $\OOO^{n+1}$ is a free $\OOO$-module, by Theorem 1 in
\cite{Frohlich75} and the consequences that follow, the ideal classes of
$\bbb$ and $\bbb'$ are equal. Up to changing $\bbb$ and $\bbb'$ in
their equivalence class, we may assume that $\bbb=\bbb'$, hence in
particular $q(z)=q(z')$.

The map $xa+yc+zb \mapsto x'a+y'c+z'b$ for all $a\in\aaa$,
$c\in\,\widecheck{\aaa}^{\;-1}$ and $b\in\bbb$ is an isomorphism of
right $\OOO$-modules from $\OOO^{n+1}$ to itself, preserving the
quaternionic Hermitian form $q$ and sending $x$ to $x'$, as wanted.
\cqfd

\bigskip
From now on, we fix an integral ideal $\mmm$ of $\OOO$, which is
bilateral and stable
%
%
by the conjugation $x\mapsto \overline{x}$, as for instance $\mmm=
(1+i)\OOO$ if $\OOO=\ZZ[\frac{1+i+j+k}{2},i,j,k]$ is the Hurwitz
maximal order in the Hamilton quaternion algebra $A= \big(\frac{-1,-1}
{\QQ}\big)$ over $\QQ$. The quotient $\ZZ$-module $\OOO/\mmm$ is then
a ring endowed with an anti-involution again denoted by $x\mapsto
\overline{x}$. We denote by $\operatorname{U}_q(\OOO/\mmm)$ the finite
group of $(n+1)\times(n+1)$ invertible matrices in $\OOO/\mmm$,
preserving the Hermitian form $-\overline{z_0}z_n- \overline{z_n}z_0 +
\sum _{i=1}^{n-1}\overline{z_i}z_i$ on $(\OOO/\mmm)^{n+1}$. Let
$\operatorname{B}_q(\OOO/\mmm)$ be its upper triangular subgroup. We
denote by $\Ga_\mmm$ the {\it Hecke congruence subgroup} of
$\operatorname{U}_q(\OOO)$ modulo $\mmm$, that is, the preimage of
$\operatorname{B}_q(\OOO/\mmm)$ by the group morphism
$\operatorname{U}_q(\OOO)\ra \operatorname{U}_q(\OOO/\mmm)$ of
reduction modulo $\mmm$. For every subgroup $H$ of
$\operatorname{U}_q$, we denote by $\operatorname{P}\!H$ its image in
$\PU_q$.

\bprop \label{prop:identiforbit} 
If $n=2$, then
\begin{enumerate}
\item 
  the set of parabolic fixed points of $\operatorname{P}\!\Ga_{\mmm}$
  is the set of points in $\partial_\infty\HH^n_\HH$, which is the
  isotropic cone of $q$ in $\PP^n_{\rm r}(\HH)$, having homogeneous
  coordinates that are elements in $\OOO$;
\item
  the orbit $\operatorname{P}\!\Ga_{\mmm}\cdot\infty$ is the set of
  points in $\partial_\infty\HH^n_\HH$ having homogeneous coordinates
  in $\PP^n_{\rm r}(\HH)$ of the form $[a:\alpha:c]$ with $(a,\alpha,
  c) \in \OOO\times\mmm^{n-1}\times\mmm$, $\tr(\overline{a}c)
  =\n(\alpha)$ and $_\OOO\langle a,\alpha,c\rangle=\OOO$;
\item
  the map which associates to $[a:\alpha:c]\in \PP^n_{\rm r}(A)$ the class
  of the left fractional ideal $_\OOO\langle a,\alpha,c\rangle$
  generated by its homogeneous coordinates induces a bijection from
  the set of cusps $\PU_q(\OOO)\bs \Par_{q,\OOO}$ of $\PU_q(\OOO)$ to
  the set of left ideal classes $_\OOO\I$ of $\OOO$.
\end{enumerate}
\eprop

The number of cusps of $\PU_q(\OOO)$ is hence exactly the class number
$h_A$ of $A$, and in particular is equal to $1$ if and only if $D_A=2,
3, 5, 7,13$. Since the simple real Lie group $\PU_q$ has rank one, the
set of ends of the quaternionic hyperbolic orbifold $\PU_q(\OOO)\bs
\HH_\HH^n$ is in bijection with the set of cusps of $\PU_q(\OOO)$, and
Theorem \ref{intro:cuspnumb} in the Introduction follows.

\medskip
\dem (1) By the previously mentioned results of Borel and
Harish-Chandra, the result is true if $\mmm=\OOO$, since any element
in $\PP^n_{\rm r}(A)$ may be represented by an element of
$\OOO^{n+1}$. As $\operatorname{P}\!\Ga_{\mmm}$ has finite index in
$\PU_q(\OOO)$, the general case follows since a discrete group and a
finite index subgroup have the same set of parabolic fixed points.

\medskip
Since $_\OOO\langle 1,0,0\rangle=\OOO$, the assertions (2) when
$\mmm=\OOO$ and (3) follows from Assertion (1) and Proposition
\ref{prop:isotropictransitive}. Assertions (2) for any $\mmm$ follows
by the definition of $\Ga_{\mmm}$, since the image of $(1,0,0)$ by a
matrix in $\GL_3(\HH)$ is its first column.  
\cqfd

\section{The covolume of $\PU_q(\OOO)$}
\label{sec:computvol}

In this section, we prove Theorem \ref{intro:computvol} in the
Introduction, using Prasad's volume formula in \cite{Prasad89} and
arguments from \cite{EmeKim18}.

Let $\P$ be the set of positive primes in $\ZZ$. For every $p\in \P$,
the order $\OOO_p=\OOO\otimes_\ZZ\ZZ_p$ is a maximal order in the
quaternion algebra $A_p=A\otimes_\QQ \QQ_p$ over $\QQ_p$ (see for
instance \cite[page 84]{Vigneras80}). Let us denote by $v_p$ the
$p$-adic valuation of $\QQ_p$ and by $\n_p$ the reduced norm on $A_p$.
For every $p\in\P$, recall that by the definition of the discriminant
$D_A$ of $A$, if $p$ does not divide $D_A$, then $A_p$ is isomorphic
to $\M_2(\QQ_p)$ and otherwise $A_p$ is the (unique up to isomorphism)
quaternion algebra over $\QQ_p$ that is a division algebra.
Furthermore, let us consider the discrete valuation $\nu_p=
\frac{1}{2}\,v_p\circ\n_p$ on $A_p$ (with value group $\frac{1}{2}
\,\ZZ$).  It coincides with $v_p$ on $\QQ_p$, which is the reason of
the factor $\frac{1}{2}$. The unique maximal order $\OOO_p$ is equal to the
valuation ring of $\nu_p$ (see for instance \cite[page
  34]{Vigneras80}). We fix a uniformiser $\pi_p\in\OOO_p$ for $\nu_p$:
we have $\n_p(\pi_p)=p$ and $\nu_p(\pi_p)=\frac{1}{2}$.

\medskip
As in the beginning of Section \ref{sec:cuspnumb}, let $\underline{G}$
be the absolutely connected, (quasi-)simple, simply connected
algebraic group over $\QQ$, endowed with a $\ZZ$-form such that
$\underline{G}(\ZZ) =\operatorname{U}_q (\OOO)$ and
$\underline{G}(\RR) =\operatorname{U}_q$. We assume that $n=2$.

Note that $\underline{G}$ is a $\QQ$-form of the split (hence
quasi-split) algebraic group $\G=\Sp_3$ over $\QQ$, whose type is
$C_{3}$, by \cite[page 89]{PlaRap94}. The $\QQ$-group $\underline{G}$
is an inner form of $\G$ since the type $C_{3}$ has no symmetries in
its diagram, by \cite[page 67]{PlaRap94}.  The {\it absolute rank} $r$
of $\G$ and the {\it exponents} $m_1,\dots, m_r$ of $\G$ are given by
\begin{equation}\label{eq:rankexpo}
r=3 \;\;\;{\rm and}\;\;\;m_1=1,\; m_2=3,\;m_3=5
\end{equation}
(see for instance \cite[page 96]{Prasad89}).

Let $\I_{\underline{G},\QQ_p}$ be the Bruhat-Tits building of
$\underline{G}$ over $\QQ_p$ (see for instance \cite{Tits79,BruTit87}
for the necessary background on Bruhat-Tits theory).  Recall that a
subgroup of $\underline{G}(\QQ_p)$ is {\it parahoric} if it is the
stabiliser of a simplex of $\I_{\underline{G},\QQ_p}$. A {\it coherent
  family of parahoric subgroups} of $\underline{G}$ is a family
$(Y_p)_{p\in\P}$, where $Y_p$ is a parahoric subgroup of
$\underline{G} (\QQ_p)$ for every $p$ and $Y_p= \underline{G} (\ZZ_p)$
for $p$ big enough. The {\it principal lattice} associated with this
family is (see \cite[\S~3.4]{Prasad89}) the subgroup of
$\underline{G}(\QQ)$ consisting of its elements which, when considered
as elements of $\underline{G}(\QQ_p)$, belong to $Y_p$, for every
$p\in\P$.

Let $p\in\P$. First assume that $p$ does not divide $D_A$. Then
$\underline{G}$ is isomorphic to $\G=\Sp_3$ over $\QQ_p$. The vertices
of the building $\I_{\G ,\QQ_p}$ are (see for instance \cite{BruTit87}
or \cite{Shemanske07}) the homothety classes of $\ZZ_p$-lattices in
${\QQ_p}^6$ generated as $\ZZ_p$-module by the union $\B$ of the
standard basis of three orthogonal hyperbolic planes, as for instance
with $\B$ the canonical basis for the standard symplectic form on
${\QQ_p}^6$.  

\blemm 
If $p$ does not divide $D_A$, then $\underline{G}(\ZZ_p)$ is parahoric.
\elemm

\dem 
In what follows, we denote by $X\mapsto\;^{t_n}\!X$ the transposition
map of $n\times n$ matrices. Note that $A_p=\M_2(\QQ_p)$ and $\OOO_p=
\M_2(\ZZ_p)$ since $p$ does not divide $D_A$, and that the conjugation
of the quaternion algebra $\M_2(\QQ_p)$ is $x= \begin{pmatrix} a &
  b\\ c & d\end{pmatrix}\mapsto x^\sigma= \begin{pmatrix} d & -b\\ -c
& a \end{pmatrix}$ (see for instance \cite[page 3]{Vigneras80}). Let
$J_0=\begin{pmatrix} 0 & -1\\ 1 & 0\end{pmatrix}$, so that for every
$x\in\M_2(\QQ_p)$, we have $x^\sigma = J_0^{-1}\;^{t_2}x\,J_0$. Let $J_1=
\begin{pmatrix} 
0 & 0 & -I_2\\ 0 & I_2 & 0\\ -I_2& 0 & 0
\end{pmatrix}$ 
and $J_2=
\begin{pmatrix} 
J_0 & 0 & 0\\ 0 & J_0 & 0\\ 0& 0 & J_0
\end{pmatrix}$. 
Considering $3\times 3$ matrices with coefficients in $\M_2(\QQ_p)$ as
$6\times 6$ matrices, an easy computation shows that for every $X$ in
$\M_3(\M_2(\QQ_p))$, with $X^\sigma$ the matrix whose coefficients are
the conjugates of the coefficients of $X$, we have $^{t_3}\!X^\sigma =
J_2^{-1}\,\;^{t_6}\!X\; J_2$. Thus $X$ belongs to $U_q(A_p)$, that is,
$^{t_3}\!X^\sigma \;J_1 \;X= J_1$, if and only if $^{t_6}\!X \;J_3
\;X= J_3$ where $J_3=J_2J_1=\begin{pmatrix} 0 & 0 & -J_0 \\ 0 & J_0 &
0\\ -J_0& 0 & 0\end{pmatrix}$.
Note that $J_3$ is (up to a harmless permutation of the canonical
basis) the matrix of the standard symplectic product defining
$\G=\Sp_3$. We hence have $\underline{G}(\ZZ_p)= \operatorname{U}_q
(\OOO_p) =\Sp_3(\ZZ_p)$, which is the stabiliser of the class of the
standard $\ZZ_p$-lattice ${\ZZ_p}^6$. The result follows. 
\cqfd

\medskip
Now assume that $p$ divides $D_A$. Then $\underline{G}(\QQ_p)=
U_q(A_p)$ has local type $^2C_{3}$ in Tits' classification, see
\cite[page 67]{Tits79}. Note that $\operatorname{SU}=\operatorname{U}$
in our case, as mentioned in \cite[Rem.2.2]{EmeKim18}.  Its
local index is shown below (see \cite[page 63]{Tits79}):
\begin{center}
\begin{picture}(0,0)%
\includegraphics{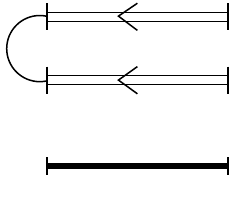}%
\end{picture}%
\setlength{\unitlength}{3812sp}%
\begingroup\makeatletter\ifx\SetFigFont\undefined%
\gdef\SetFigFont#1#2#3#4#5{%
  \reset@font\fontsize{#1}{#2pt}%
  \fontfamily{#3}\fontseries{#4}\fontshape{#5}%
  \selectfont}%
\fi\endgroup%
\begin{picture}(1166,1066)(848,-395)
\put(863,-215){\makebox(0,0)[lb]{\smash{{\SetFigFont{11}{13.2}{\rmdefault}{\mddefault}{\updefault}{\color[rgb]{0,0,0}$\times$}%
}}}}
\put(1036,-331){\makebox(0,0)[lb]{\smash{{\SetFigFont{11}{13.2}{\rmdefault}{\mddefault}{\updefault}{\color[rgb]{0,0,0}$s$}%
}}}}
\put(1032,-61){\makebox(0,0)[lb]{\smash{{\SetFigFont{11}{13.2}{\rmdefault}{\mddefault}{\updefault}{\color[rgb]{0,0,0}$3$}%
}}}}
\put(1936,-61){\makebox(0,0)[lb]{\smash{{\SetFigFont{11}{13.2}{\rmdefault}{\mddefault}{\updefault}{\color[rgb]{0,0,0}$2$}%
}}}}
\put(1936,-331){\makebox(0,0)[lb]{\smash{{\SetFigFont{11}{13.2}{\rmdefault}{\mddefault}{\updefault}{\color[rgb]{0,0,0}$s$}%
}}}}
\end{picture}%

\end{center}
In particular, $\underline{G}(\QQ_p)$ has relative rank $1$. By
\cite[Ex.~2.7]{Tits79}, the building $\I_{\underline{G}, \QQ_p}$ is a
biregular tree of degrees $p^3+1$ and $p^2+1$, and all its vertices
are special.

\blemm 
If $p$ divides $D_A$ and $p\neq 2$, then $\underline{G}(\ZZ_p)$ 
is parahoric.
\elemm

See also \cite[Lem.~5.6]{EmeKim18} with a different proof, our proof
being useful in order to deal with the case $p=2\;|\,D_A$.

\medskip\dem 
We will use the interpretation of $\I_{\underline{G}, \QQ_p}$ as the
set of selfdual norms on $A_p^{\,3}$ (see \cite{BruTit84, BruTit87},
which uses a $-\operatorname{log}_p$ version of them, in order to
allow for infinite residual fields).

With the notation of \cite[\S 1.2]{BruTit87}, we take $K=L:=\QQ_p$,
$D:=A_p$, $x\mapsto x^\sigma:=\ov x$ the quaternion conjugation in
$A_p$, $\varepsilon:=+1$ so that $D^0=\{x\in A_p\;:\; \tr\;x=0\}$,
$X:=A_p^{\,3}$ considered as a right vector space over $A_p$, and
$b:X\times X\ra D$ the form induced by extension of scalars to $\QQ_p$
of the restriction $\Phi_A:A^3\times A^3\ra A$ to $A$ of the quaternionic
sesquilinear form $\Phi$ defined in Equation \eqref{eq:Phi} with
$n=2$,  so that the form $q:X\ra D/D^0=K$ of loc.~cit.~coincides with
the extension of scalars to $\QQ_p$ of our form $q:A^3\ra \QQ$ over
$\QQ$ see Equation \eqref{eq:q}.

With the notation of \cite[\S 1.2]{BruTit87}, we take $\omega_L=
\omega:=v_p$, which is a discrete normalised valuation on $K=L=\QQ_p$
and $\omega_D=\nu_p$, which is a discrete valuation (with value group
$\frac12\ZZ$) on $D=A_p$ extending $v_p$.

Let $\lambda\mapsto |\lambda|_p=p^{-\nu_p(\lambda)}$, which is a map
from $A_p$ to $[0,+\infty[$, be the unique extension to $A_p$ of the
absolute value of $\QQ_p$. A {\it norm}\footnote{Compare with \cite[\S
1.1]{BruTit84}: a map $\alpha:V\ra [0,+\infty[$ is a norm as defined 
here if and only if $-\operatorname{log}_p \alpha$ is a norm in the
sense of \cite{BruTit84}.}  on the right $A_p$-vector space
$V=A_p^{\,3}$ is a map $\alpha:V\ra [0,+\infty[$ such that

$\bullet$~ $\alpha(x)=0$ if and only if $x=0$

$\bullet$~ $\alpha(x\lambda)=|\lambda|_p\; \alpha(x)$ for all $x\in V$
        and $\lambda\in A_p$,

$\bullet$~ $\alpha(x+y)\leq \max\{\;\alpha(x),\;\alpha(y)\,\}$ for all
        $x,y\in V$.

\noindent 
As defined in \cite[\S 2.5]{BruTit87}, its {\it dual norm} is the norm
$$
\overline{\alpha}:x\mapsto \max_{y\in V-\{0\}}\;\frac{|\,f(x,y)\,|_p}
{\alpha(y)}
$$ 
on $V$, where $f$ (keeping the notation of \cite[\S 1.2]{BruTit87}) is
the map from $V\times V$ to $\QQ_p$ defined by $f(x,y)=q(x+y)-q(x) -
q(y)$. A norm $\alpha$ on $V$ is {\it selfdual} if $\alpha=
\overline{\alpha}$. The linear action of $\underline{G}(\QQ_p)=
\operatorname{U}_q(A_p)$ on $V$ induces a left action on the set of
selfdual norms on $V$ by $(g,\alpha)\mapsto \alpha\circ g^{-1}$.

A {\it Witt basis} of $V$ is a basis $(e_{-1},e_0,e_{1})$ of the right
$A_p$-vector space $V$ such that 
$$
q(e_{\pm 1})=f(e_0,e_{\pm 1})=0,\;\;\; {\rm and} \;\;\;
q(e_0)=f(e_1,e_{-1})=1\;.
$$ 
By \cite[\S 2]{BruTit87},\footnote{More precicely, this follows from
  Theorem 2.12, 2.5 Proposition (ii) and 2.9 Proposition of op.~cit.,
  using the fact that the value group of $\nu_p$ is $\frac{1}{2}\ZZ$.}
there exists a $\underline{G}(\QQ_p)$-equivariant bijection from the
building $\I_{\underline{G}, \QQ_p}$ to the set of selfdual norms on
$V$ such that, for every Witt basis $(e_{-1},e_0,e_{1})$ of $V$, the
sequence $(\alpha_n)_{n\in\ZZ}$ of norms
$$
\alpha_n:x=\sum_{i=-1}^1 e_i\lambda_i\mapsto  
\max\{\;p^{-\frac{n}{2}}|\lambda_{-1}|_p,\;|\lambda_{0}|_p,\;
p^{\frac{n}{2}}|\lambda_{1}|_p\;\}\,,
$$ 
for $n\in\ZZ$, is the sequence of selfdual norms associated with the
sequence of vertices $(x_n)_{n\in\ZZ}$ along an appartment of
$\I_{\underline{G}, \QQ_p}$, such that $\alpha_0$ is associated with
$x_0$.  Furthermore, let
$$
{\cal X}_n=e_{-1}\pi_p^{-n}\OOO_p+ e_0\OOO_p +e_{1}\pi_p^{n}\OOO_p
$$ 
be the right $\OOO_p$-lattice generated by the Witt basis
$(e_{-1}\pi_p^{-n}, e_0\OOO_p,e_{1}\pi_p^{n})$, which is the unit ball
of the norm $\alpha_n$, since $|\pi_p|_p=p^{-\frac12}$. Then by \S 3.9
page 180 of \cite{BruTit87}, the smooth affine group scheme $\G_{x_0}$
over $\ZZ_p$ associated with the vertex $x_0$ is the schematic closure
of $\underline{G}$ in the $\ZZ_p$-form of the general linear group
$\GL_3(A_p)$ over $\QQ_p$ defined by the $\ZZ_p$-lattice ${\cal X}_0$,
and $\G_{x_0}(\ZZ_p)$ is the stabiliser of the vertex $x_0$ in
$\underline{G}(\QQ_p)$.

Since $p\neq 2$, the element $-2$ is an invertible element of $\ZZ_p$
hence of $\OOO_p$. If $(e'_{-1},e'_0,e'_{1})$ is the canonical basis
of $V=A_p^{\,3}$, which satisfies $q(e'_{\pm 1})=f(e'_0,e'_{\pm 1})
=0$, $q(e'_0)=1$ and $f(e'_1,e'_{-1})=-2$, then $(e_{-1},e_0,e_{1}) =
(e'_{-1},e'_0,e'_{1}\,\frac{1}{-2})$ is a Witt basis of $V$, and
generates the same right $\OOO_p$-lattice ${\cal X}_0$ as the
canonical basis. Therefore, $U_q(\OOO_p)=\underline{G}(\QQ_p)
\cap\GL_3(\OOO_p) =\G_{x_0}(\ZZ_p)$ is parahoric. 
\cqfd

\brema\label{rem:caspdeux}
When $p$ divides $D_A$ and $p=2$, the group $U_q(\OOO_p)$ is not
parahoric.
%
%
Indeed, again with $(e'_{-1},e'_0,e'_{1})$ the canonical basis of $V$,
since $\n_2(\pi_2)=2$, the basis $(e_{-1},e_0,e_{1}) = (e'_{-1}
\pi_2^{-1},e'_0,-e'_{1}\pi_2^{-1})$ is a Witt basis of
$V$. With the above notation, $U_q(\OOO_2)$ is the subgroup of the
stabiliser of $x_0$ in $\underline{G}(\QQ_2)$ fixing the two edges
with origin $x_0$ and endpoints $x_{\pm 1}$, since ${\cal X}_{-1}\cap
{\cal X}_{1}= e'_{-1}\OOO_p+e'_0\,\OOO_p+e'_{1}\,\OOO_p$.
\erema

Thus, by definition, if $D_A$ is odd, the family
$(\underline{G}(\ZZ_p))_{p\in\P}$ is a coherent family of (maximal)
parahoric subgroups of $\underline{G}$, and $\operatorname{U}_q
(\OOO)= \underline{G}(\ZZ)= \{g\in \underline{G}(\QQ)\;:\;\forall
p\in\P, g\in \underline{G}(\ZZ_p)\}$ is its associated principal
lattice. If $D_A$ is even, the family $(Y_p)_{p\in\P}$ with
$Y_p=\underline{G}(\ZZ_p)$ if $p\neq 2$ and $Y_2$ the stabiliser in
$\underline{G}(\QQ_2)$ of the point $x_0$ defined in Remark
\ref{rem:caspdeux}, is a coherent family of (maximal) parahoric
subgroups of $\underline{G}$. We will compute below the index of
$\operatorname{U}_q (\OOO)$ in the associated principal lattice when
$D_A$ is even.

\medskip 
For every $p\in\P$,

$\bullet$~~ let $y_p$ (respectively ${\mathfrak y}_p$) be the vertex
of $\I_{\underline{G}, \QQ_p}$ (respectively $\I_{\G,\QQ_p}$)
stabilised by the subgroup $\operatorname{U}_q(\OOO_p)$
(respectively $\Sp_3(\ZZ_p)$), such that if $D_A$ is even, then $y_2$
is the point $x_0$ defined in Remark \ref{rem:caspdeux}.

$\bullet$~~ let $\ov M_p$ (respectively $\ov\M_p$) be the maximal
reductive quotient, defined over the residual field $\FF_p=
\ZZ_p/p\ZZ_p$, of the identity component of the reduction modulo $p$
of the smooth affine group scheme over $\ZZ_p$ associated with $y_p$
(respectively ${\mathfrak y}_p$); see for instance \cite[\S
  3.5]{Tits79} with $\Omega=\{v\}$.
  
Note that $\ov M_p=\ov\M_p$ if $p$ does not divide
$D_A$, and that for every $p\in\P$ the algebraic group $\ov\M_p$ is
isomorphic to $\Sp_3$ (of type $C_3$) over $\FF_p$. In particular
$\ov\M_p(\FF_p)=\Sp_3(\FF_p)$ and thus, for every $p\in\P$, the orders
of the finite groups of Lie type being listed for example
in~\cite[Table~1]{Ono66}, we have
\begin{equation}\label{eq:dimsplit}
\dim\; \ov\M_p=21 \;\;{\rm and}\;\;
|\;\ov \M_p(\FF_p)\,|=p^9(p^2-1)(p^4-1)(p^6-1)\;.
\end{equation}

Assume now that $p$ divides $D_A$ and $p\neq 2$.  Note that the pair
$(L=\OOO^3,h=\Phi)$, where $\Phi$ is defined by Equation
\eqref{eq:Phi} with $n=2$, is {\it admissible of maximal type} over
$p$ in the sense of \cite[\S 5.1 and Def.~5.3]{EmeKim18}.  Indeed,
$\OOO$ is a maximal order in $A=\OOO\otimes_\ZZ\QQ$. The pair $(L,h)$
is a Hermitian right $\OOO$-module with Witt signature $(1,2)$, which
is unimodular (called {\em regular} in \cite{EmeKim18}) over $p$.
Recall that this means that the map from $\widecheck{\OOO_p^{\,3}}$ to
$\operatorname{Hom} _{\OOO_p} (\OOO_p^{\,3},\OOO_p)$ is an isomorphism
of left $\OOO_p$-modules.  This property is indeed satisfied by
\cite[Lem.~5.1]{EmeKim18} since $p\neq 2$. This restriction $p\neq 2$
is needed since putting $h$ in diagonal form as in
\cite[Eq.~(5.1)]{EmeKim18} requires to invert $2$.

By \cite[Lem.~4.1]{EmeKim18} with $k=\QQ$, $v=p$, ${\bf G}
=\underline{G}$, $t=0$ and $P^t_v=\operatorname{U}_q(\OOO_p)$, whose
hypothesis is satisfied by \cite[Lem.~5.6]{EmeKim18}, we have
\begin{equation}\label{eq:dimnonsplit}
p^{(\dim \,\ov M_p-\dim \,\ov \M_p)/2}\;
\frac{|\,\ov \M_p(\FF_p)|}{|\,\ov M_p(\FF_p)|}
= (p-1)(p^2+1)(p^3-1)\;.
\end{equation}

Assume finally that $p=2$ divides $D_A$. The Tits index of $\ov \M_p$, as
computed by the rule of \cite[3.5.2]{Tits79} is $^2A_2$, and by
\cite[3.5.4]{Tits79}, the link of the vertex $y_p$ in the tree
$\I_{\underline{G}, \QQ_p}$ canonically identifies with the spherical
building of $\ov \M_p$ over $\FF_p$. By \cite[page 55]{Tits66}, $\ov
\M_p$ is hence
the group $\operatorname{U}_q(\FF_{p^2})$ where the involution on
$\FF_{p^2}$ is its Frobenius automorphism $x\mapsto \overline{x}=
x^p$. The spherical building of $\ov \M_p$ is the finite set of
isotropic points in the projective plane over $\FF_{p^2}$. The vertices of
the link of $y_2=x_0$ corresponding to $x_{-1}$ and $x_1$ with the
notation of Remark \ref{rem:caspdeux} are the projective points
defined by the (isotropic) first and last vectors of the canonical
basis of $(\FF_{p^2})^3$. Let $H$ be the intersection of the
stabilisers in $\operatorname{U}_q(\FF_{p^2})$ of the two isotropic
points $[1:0:0]$ and $[0:0:1]$. An easy computation shows that $H$
consists of the diagonal matrices $\begin{pmatrix} a & 0 & 0\\ 0 & U &
0\\0 & 0 & d\end{pmatrix}$, with $a,U,d\in \FF_{p^2}^{\;\times}$,
$d=\frac{1}{\overline{a}}=a^{-p}$ and $U^{p+1}=U\,\overline{U} =1$.
Since the multiplicative group $\FF_{p^2}^\times$ is isomorphic to the
additive cyclic group $\ZZ/((p^2-1)\ZZ)$, which contains exactly $p+1$
elements $x$ such that $(p+1)x=0$, the order of $H$ is equal to
$(p^2-1)(p+1)$. The center of $\operatorname{U}_q(\FF_{p^2})$ has
order $p+1$, and the quotient by its center is the finite simple group
called the {\it Steinberg group} $^2A_2(p^2)$, whose order is
$p^3(p^2-1)(p^3+1)$. Hence the index of $H$ in
$\operatorname{U}_q(\FF_{p^2})$ is $p^3(p^3+1)$. By Remark
\ref{rem:caspdeux} and since $p=2$, we hence have that the index of
$U_q(\OOO_2)$ in the parahoric subgroup $Y_2$ is
\begin{equation}\label{eq:indexchiant}
[Y_2:U_q(\OOO_2)]=[U_q(\FF_{2^2}):H]=72\;.
\end{equation}

\medskip
Now, let $\mu$ be the Haar measure on $\underline{G}(\RR)=
\operatorname{U}_q=\Sp(1,2)$ normalized as in \cite[\S 3.6]{Prasad89}.
The next lemma relates it to the Riemannian measure coming from the
choice made in Section \ref{sec:quathypspace} of the sectional
curvature on $\HH^2_\HH$.

\blemm \label{lem:tamagawa} We have
$
\Vol(\PU_q(\OOO_q)\,\bs \,\HH^2_\HH)=\frac{\pi^4}{120}\;
\mu(\operatorname{U}_q(\OOO_q)\,\bs\operatorname{U}_q)\;.
$
\elemm

\dem The proof is similar to the one in \cite{Emery09} or Emery's
appendix of \cite{ParPau13ANT}.

By the definition of $\mu$, if $w$ is the top degree exterior form on
the real Lie algebra of $\underline{G}(\RR)$ whose associated
invariant differential form on $\underline{G}(\RR)$ defines the
measure $\mu$ and if $G_u=\underline{G}_u(\RR)$ is the compact real form
of $\underline{G}(\CC)$, then the complexification $w_\CC$ of $w$ on
the complex Lie algebra of $\underline{G}(\CC)=\underline{G}_u(\CC)$
defines a top degree exterior form $w_u$ on the real Lie algebra of
$G_u$, whose associated invariant differential form on $G_u$ defines a
measure $\mu_u$, and we require that $\mu_u(G_u)=1$.

Let $\mu'$ be the Haar measure on the noncompact real Lie group
$\Sp(1,2)$ that disintegrates by the fibration $\Sp(1,2)\ra
\Sp(1,2)/(\Sp(1)\times\Sp(2))=\HH^2_\HH$ with measures on the fibers
of total mass one and measure on the base the Riemannian measure
$d{\rm vol}_{\HH^2_\HH}$ of the Riemannian metric with sectional
curvatures contained in $[-4,-1]$, as in Section
\ref{sec:quathypspace}. In particular,
$$
\Vol(\PU_q(\OOO_q)\,\bs\, \HH^2_\HH)=
\mu'(\operatorname{U}_q(\OOO_q)\,\bs\operatorname{U}_q)\;.
$$

Let $\mu'_u$ be the Haar measure on the compact real Lie group
$\Sp(3)$ that disintegrates by the fibration $\Sp(3)\ra
\Sp(3)/(\Sp(1)\times\Sp(2))=\PP^2_{\rm r}(\HH)$ with measures on the
fibers of total mass one and measure on the base the Riemannian
measure $d{\rm vol}_{\PP^2_{\rm r}(\HH)}$ of the Riemannian metric
with sectional curvatures contained in $[1,4]$. By
\cite[Ex.A.III.8]{BerGauMaz71}, this Riemannian metric is the standard
Fubini-Study metric, and
$$
\mu'_u(\Sp(3))=\Vol(\PP^2_{\rm r}(\HH))=\frac{\pi^4}{120}\;.
$$ 

The duality between irreducible symmetric spaces of noncompact type
endowed with a left invariant Riemannian metric and the ones of
compact type sends $\HH^2_\HH$ to $\PP^2_{\rm r}(\HH)$, with opposite
signs on the range of the sectional curvatures (see for instance
\cite[Ch.~5]{Helgason78}), and hence $\mu' = \frac{\pi^4}{120}
\,\mu$. The result follows.  \cqfd

\bigskip
We now apply Prasad's volume formula \cite[Theo.~3.7]{Prasad89}. With
the notation of this theorem, we take 

$\bullet$~ $k=\QQ$ 

\noindent so that its set of infinite places is $V_\infty=
\{\infty\}$, its set of finite places is $V_f=\P$, and the order $q_v$
of the residual field $\fff_v=\ZZ_p/p\ZZ_p=\FF_p$ is $p$ for every
$v=p\in V_f$,

$\bullet$~ $S= \{\infty\}$ 

\noindent so that $S_f=S\cap V_f$ is empty and

$\bullet$~ ${\bf G}=\underline{G}$, 

\noindent so that the Tamagawa number $\tau_k({\bf G})$ is $1$ by
page 109 of op.~cit.~since $k$ is a number field, $\ell=\QQ$ is a
smallest splitting field of $\G$ over $\QQ$ (since $\G$ is split over
$\QQ$), and the discriminants of $k$ and $\ell$ over $\QQ$ are
$D_k=D_\ell=1$. We hence have, since $\ov\M_p=\ov M_p$ if $p$ does not
divide $D_A$ and by Equation \eqref{eq:rankexpo} for the second
equality,
\begin{align}
\mu(\operatorname{U}_q(\OOO)\bs\operatorname{U}_q)&=\prod_{i=1}^{r}
\frac{(m_i)!}{(2\pi)^{m_i+1}}\prod_{p\in\P}
\frac{p^{(\dim \,\ov M_p+\dim
    \,\ov \M_p)/2}}{|\,\ov M_p(\FF_p)|}\nonumber\\ 
\label{eq:formulprasad}&=
\frac{720}{(2\pi)^{12}}\prod_{p\in\P}\frac{p^{\dim \,\ov \M_p}}{|\ov
    \M_p(\FF_p)|}\;
\prod_{p|D_A}\frac{|\,\ov \M_p(\FF_p)|}{|\,\ov M_p(\FF_p)|}\;
p^{(\dim \ov M_p-\dim \ov \M_p)/2}\;.
\end{align}
Using Euler's product formula $\zeta(s)= \prod_{p\in\P}
\frac{1}{1-p^{-s}}$ for Riemann's zeta function, we have by Equation
\eqref{eq:dimsplit}, since $\zeta(2)=\frac{\pi^2}{6}$,
$\zeta(4)=\frac{\pi^4}{90}$ and $\zeta(6)=\frac{\pi^6}{945}$,
\begin{equation}\label{eq:prodzeta}  
\prod_{p\in\P}\frac{p^{\dim \,\ov \M_p}}{|\,\ov
    \M_p(\FF_p)|}=\zeta(2)\,\zeta(4)\,\zeta(6)=\frac{\pi^{12}}{510300}\;.
\end{equation}

Theorem \ref{intro:computvol} in the Introduction when $D_A$ is odd
follows from Equations \eqref{eq:formulprasad}, \eqref{eq:prodzeta},
\eqref{eq:dimnonsplit}, and from Lemma \ref{lem:tamagawa}. When $D_A$
is even, the obtained formula computes the covolume of the principal
lattice $\Lambda$ associated with the coherent family
$(Y_p)_{p\in\P}$. By construction, $\operatorname{U}_q(\OOO)$ is
exactly the subgroup of elements in $\Lambda$ which, when considered
in $\underline{G}(\QQ_2)$, belong to the finite index subgroup
$\operatorname{U}_q(\OOO_2)$ of $Y_2$. Since $\Lambda$ is dense in
$Y_2$ and $\operatorname{U}_q(\OOO)$ is dense in
$\operatorname{U}_q(\OOO_2)$, this proves that the index $m_A$ of
$\operatorname{U}_q(\OOO)$ in $\Lambda$ is equal to the index of
$\operatorname{U}_q(\OOO_2)$ in $Y_2$, which is $72$ by Equation
\eqref{eq:indexchiant} if $D_A$ is even.


\section{Horospherical quaternionic hyperbolic geometry}
\label{sec:quathypgeom}

In this section, we describe the geometry of the horospheres in the
quaternionic hyperbolic space $\hnh$ (see also \cite{KimPar03,
  Philippe16}). We introduce the quaternionic Heisenberg group and
discuss the geometry of its quaternionic contact structure (see for
instance \cite{Biquard01}).

The {\it horospherical coordinates} $(\zeta,u,t)\in \HH^{n-1} \times
\Im\;\HH \times \, [0,+\infty[$, that we will use from now on unless
otherwise stated, of $(w_0,w)\in\hnh\cup(\partial_\infty \hnh-
\{\infty\})$ are
%
%
\begin{equation}
\begin{aligned}\label{eq:horosphecoord}
(\zeta,u,t) & =(w,\;2\,\Im\; w_0, \;\tr\, w_0-\n(w))\\
{\rm hence}\;\;\;
(w_0,w) & =\Big(\;\frac{\n(\zeta)+t+u}{2},\;\zeta\,\Big)\,,
\end{aligned}
\end{equation}
so that the Riemannian metric of $\hnh$ is given by 
\begin{equation}\label{eq:riemethorosphecoord}
ds^2_{\,\hnh}= \frac{1}{4\,t^2}\big(dt^2
+\n(du-2\,\Im \;\overline{d\zeta}\cdot\zeta)+
4\,t\,\n(d\zeta) \,\big)\,.
\end{equation} 
In horospherical coordinates, the geodesic lines from
$(\zeta,u,0)\in\partial_\infty \hnh- \{\infty\} $ to $\infty$ are, up
to translations at the source, the maps $s\mapsto (\zeta,u,e^{2s})$,
by the normalisation of the metric.

The {\em Busemann cocycle}\index{Busemann
  cocycle}\index{cocycle!Busemann} of $\hnh$ is the map
$\beta: \partial_{\infty} \hnh\times\hnh\times\hnh\to\RR$ defined by
$$
(\xi,x,y)\mapsto \beta_{\xi}(x,y)=
\lim_{s\to+\infty}d(\xi_s,x)-d(\xi_s,y)\;,
$$
where $s\mapsto \xi_s$ is any geodesic ray ending at $\xi$. It is
invariant under the diagonal action of the isometry group of $\hnh$.
The {\em horosphere} with centre $\xi \in \partial_{\infty}\hnh$
through $x\in \hnh$ is $\{y\in \hnh\;:\; \beta_{\xi}(x,y)=0\}$, and
$\{y\in \hnh\;:\; \beta_{\xi}(x,y)\leq 0\}$ is the (closed) {\em
  horoball} centred at $\xi$ bounded by this horosphere.

Given two points $x=(\zeta,u,t)$ and $x'=(\zeta',u',t')$ in $\hnh$,
the maps $\xi_s: s\mapsto (\zeta,u,e^{2s})$ and $\xi'_s:s\mapsto
(\zeta',u',e^{2s})$ are geodesic lines in $\hnh$ through $x$ and $x'$
respectively, converging to $\infty$ as $s\ra+\infty$. The Riemannian
length of the affine path from $\xi_s$ to $\xi'_s$ is bounded by a
constant times $e^{-s}$, hence $\lim_{s\ra+\infty}
d(\xi_s,\xi'_s)=0$. Thus
\begin{equation}\label{eq:busemaninfini}
\beta_\infty(x,x')=\frac{1}{2}\,\ln \frac{t'}{t}\;.
\end{equation}
The closed horoballs centred at $\infty\in\partial_\infty \hnh$ are
therefore the subsets
\begin{equation}\label{eq:defhoroinfty}
\H_s=\{(\zeta,u,t)\in\hnh\;:\;t\geq s\}= 
\{(w_0,w)\in\hnh\;:\;\tr w_0-\n(w)\geq s\},
\end{equation}
and the horospheres centred at $\infty$ are their boundaries
$\partial\H_s$, where $s$ ranges in $]0,+\infty[\,$.  Note that, for
    every $s\geq 1$, we have
\begin{equation}\label{eq:distentrhorob}
d(\partial\H_1,\partial\H_s)=\frac{\ln s}{2}\;.
\end{equation}

The {\it Cygan distance}\footnote{It is analogous to the Euclidean
  distance on the closure in $\RR^n$ of the upper halfspace model of
  $\hnr$.} on $\hnh\cup (\partial_\infty \hnh-\{\infty\})$ is (see for
instance \cite{KimPar03})
\begin{equation}\label{eq:defidistCyg}
d_{\rm Cyg}((\zeta,u,t),(\zeta',u',t'))=
\n\big(\,\n(\zeta-\zeta')+|t-t'|+
(u-u'-2\,\Im\;\overline{\zeta}\,\cdot\zeta')\;\big)^{1/4}\,.
\end{equation}

The {\it quaternionic Heisenberg group} $\HHeis_{4n-1}$ of dimension
$4n-1$ is the real Lie group structure on $\HH^{n-1}\times \Im\;\HH$
with law
$$
(\zeta,u)(\zeta',u')=
(\zeta+\zeta',u+u'+2\,\Im\;\overline{\zeta}\cdot\zeta')
$$ 
and inverses $(\zeta,u)^{-1}=(-\zeta,-u)$. When $n=2$, using the
change of coordinates $\zeta =w$ and $u=2\,\Im\,w_0$ as in Equation
\eqref{eq:horosphecoord} with $t=0$, we recover the definition given
in the Introduction. The group $\HHeis_{4n-1}$ identifies with
$\partial_\infty \hnh- \{\infty\}$ by the map $(\zeta,u)\mapsto
(\zeta,u,0)$. It identifies with a subgroup of $\operatorname{PB}_q
\subset\PU_q$ by $(\zeta,u)\mapsto \pm\begin{pmatrix} 1 & \zeta^* &
\frac{\n(\zeta)+u}{2}\\0 & I_{n-1} & \zeta\\0 & 0 & 1 \end{pmatrix}$,
where $\zeta\in\M_{n-1,\,1}(\HH)$ also denotes the column vector of
$\zeta\in\HH^{n-1}$. It acts on the space $\hnh\cup (\partial_\infty
\hnh- \{\infty\})$ by the {\it Heisenberg translations}
$$
(\zeta,u)(\zeta',u',t')= (\zeta+\zeta',\;u+u'+ 2\,\Im\;
\overline{\zeta}\cdot\zeta',\;t')\;.
$$
They are isometries for both the Riemannian metric and the Cygan
distance, and they preserve the horospheres centred at $\infty$.  For
every $u\in\Im\;\HH$, the Heisenberg translation by $(0,u)$ is called a
{\it vertical translation}.

It is easy to see that the Cygan distance on $\HHeis_{4n-1}$ is the
unique left-invariant distance on $\HHeis_{4n-1}$ with
$$
d_{\rm Cyg} ((\zeta,u),(0,0))= (\n(\zeta)^2+\n(u))^{\frac{1}{4}}\;,
$$
or equivalently using Equation \eqref{eq:horosphecoord} that if
$(w_0,w)\in\partial_\infty \hnh- \{\infty\}$, then
$$
d_{\rm Cyg} ((w_0,w),(0,0))= (4\n(w_0))^{\frac{1}{4}}\;.
$$

\medskip 
We conclude this section with geometric lemmas that will be useful in
Sections \ref{sect:measurecomput} and \ref{sect:mertens}. See also
\cite[\S 3]{Kim15}, with slightly different conventions, for a
computation similar to Lemma \ref{lem:busemann}.  The proofs are
analogous to those in the complex hyperbolic case with the added
ingredient of being careful with the noncommutativity of
multiplication in the quaternionic case.

\blemm\label{lem:busemann}
For all $x=(\zeta,u,t)$ and $x'=(\zeta',u',t')$ in $\hnh$, for all
$(\xi,r)\in \HHeis_{4n-1}=\partial_\infty \hnh-\{\infty\}$, we have
$$
\beta_{(\xi,\,r)}(x,x')=\frac{1}{2}\,\ln 
\frac{t'\;d_{\rm Cyg}(x,(\xi,r))^4}{t\;d_{\rm Cyg}(x',(\xi,r))^4}\,.
$$
\elemm

\dem 
It is easy to check that the map $\iota:(w_0,w)\mapsto (w_0^{\;-1},
ww_0^{\;-1})$ is an isometric involution of $\hnh$ sending $(0,0) \in
\partial_\infty\hnh$ to $\infty$, induced by $\begin{pmatrix} 0 & 0 &
  1\\0 & I_{n-1} & 0\\1 & 0 & 0 \end{pmatrix}$, which does belong to
$\operatorname{U}_q$.  Hence, with $x=(w_0,w)$ and $x'=(w'_0,w')$,
using Equations \eqref{eq:busemaninfini} and \eqref{eq:horosphecoord}
and the fact that $d_{\rm Cyg}(x,(0,0))^4=4\n(w_0)$ and $d_{\rm Cyg}
(x',(0,0))^4 = 4\n(w_0')$, we have
\begin{align*}
\beta_{(0,\,0)}(x,x')&=\beta_{\iota(0,\,0)}(\iota x,\iota x')=
\frac{1}{2}\,\ln \frac{\tr (w'_0)^{-1}-\n(w'(w'_0)^{-1})}
{\tr w_0^{\;-1}-\n(ww_0^{\;-1})}
\\ &= \frac{1}{2}\,\ln \frac{t'\n((w'_0)^{-1})}
{t\n(w_0^{\;-1})} =\frac{1}{2}\,
\ln\;\frac{t'\;d_{\rm Cyg}(x,(0,0))^4}{t\;d_{\rm Cyg}(x',(0,0))^4}\;.
\end{align*}
The Heisenberg translation $\tau$ by $(\xi,r)$ preserves the last
horospherical coordinates and the Cygan distances. We have
$\beta_{(\xi,\,r)}(x,x')= \beta_{(0,\,0)} (\tau^{-1}x, \tau^{-1}x')$,
since $\tau$ is an isometry of $\hnh$. This proves Lemma
\ref{lem:busemann}.  
\cqfd

\blemm \label{lem:orthprojgeod}
The orthogonal projection from $\partial_\infty\hnh-\{(0,0), \infty\}$
to the geodesic line in $\hnh$ with points at infinity $(0,0)$ and
$\infty$ is $(w_0,w)\mapsto (2\,\n(w_0)^{\frac{1}{2}},\;0)$, that is,
in horospherical coordinates, $(\zeta,u,0)\mapsto \big(0,0,
(\,\n(\zeta)^2+\n(u)\,)^{1/2}\big)$.  
\elemm

In particular, the preimages by this orthogonal projection are the
spheres of center $(0,0)$ for the Cygan distance on $\HHeis_{4n-1}$.

\medskip \dem 
For every parameter $a$ ranging in $]0,+\infty[\,$, consider the
horosphere $\partial\H_a$ centred at $\infty$. Its image by the
isometric involution $\iota:(w_0,w)\mapsto (w_0^{\;-1}, ww_0^{\;-1})$
is, using Equation \eqref{eq:horosphecoord}, the horosphere
$\{\,(\xi,r,t)\in\hnh\;:\;t=\frac{a}4\big((\n(\xi)+t)^2+ \n(r)\big)\,\}$
centred at $(0,0)$. The image of this horosphere by the Heisenberg
translation by $(\zeta,u)$ is the horosphere
$$
\{\,(\xi,r,t)\in\hnh\;:\;t=\frac{a}4\,
\big((\n(\xi-\zeta)+t)^2+\n(r-u-2\,\Im
\;\overline{\zeta}\cdot\xi\big))\,\}
$$ 
centred at $(\zeta,u)$. The orthogonal projection of $(\zeta,u)$ on
the geodesic line $\ell$ from $(0,0)$ to $\infty$ is attained when the
parameter $a$ gives a double point of intersection $(0,0,t)$ between
this horosphere and $\ell$. The quadratic equation 
$$
t=\frac{a}4\,((\n(\zeta)+t)^2+ \n(u))
$$ 
whose unknown is $t$ has a double solution if and only if its reduced
discriminant $\Delta'= (\n(\zeta)-\frac{2}{a})^2- (\n(\zeta)^2
+\n(u))$ vanishes, that is, since $a>0$, if and only if $a=
\frac{2}{(\n(\zeta)^2 +\n(u))^{1/2} + \n(\zeta)}$, giving
$t=(\n(\zeta)^2 +\n(u))^{1/2}$.  The result follows. \cqfd

\blemm \label{lem:projsurquatgeodline} 
Let $C$ be the quaternionic geodesic line
$\{(w_0,w)\in\hnh\;:\;w=0\}$. The orthogonal projection from $\hnh$ to
$C$ is the map $(w_0,w)\mapsto (w_0,0)$. On $\partial_\infty \hnh-
\partial_\infty C$ endowed with the horospherical coordinates, this
map extends as $(\zeta,u,0)\mapsto (0,u,\n(\zeta))$.  
\elemm

\dem 
Let $(w_0,w)\in\hnh$. It is easy to check that the distance
given (see Equation \eqref{ed:distexplicit}) by the formula
$$
\cosh^2d((w_0,w),(w'_0,0))=
\frac{\n(\,\overline{w_0}+w'_0)}{-q(w_0,w,1)\;\tr w'_0}
$$ 
is minimised over  $(w'_0,0)\in C$ exactly when $w'_0=w_0$. 

Since $C$ is totally geodesic, the closest point mapping from $\hnh$
to $C$ coincides with the orthogonal projection, which is hence
$(w_0,w)\mapsto (w_0,0)$. The expression in horospherical coordinates
of the boundary extension follows from the equations in
\eqref{eq:horosphecoord}.  
\cqfd

\blemm\label{lem:calcgeodline} For every $(w_0,w)\in\partial_\infty
\hnh-\{\infty\}$ with $w_0\neq 0$, the map from $\RR$ to $\hnh$
defined by $s\mapsto (w_0(1+2e^{2s}w_0)^{-1},\;w(1+2e^{2s}w_0)^{-1})
\in \hnh$ is a geodesic line from $(w_0,w)$ to $(0,0)$.  
\elemm

\dem 
The image of $\begin{pmatrix} 1 & 0 & 0\\ ww_0^{-1} & I_{n-1} &
  0\\ w_0^{-1} & (ww_0^{-1})^* & 1 \end{pmatrix}$ in
$\operatorname{PGL}_{n+1}(\HH)$ is the conjugate by the isometric
involution $\iota:(w'_0,w')\mapsto ({w'_0}^{\;-1}, w'{w'_0}^{\;-1})$
of the Heisenberg translation by $\iota(w_0,w)$. Hence it belongs to
$\PU_q$, fixes $\iota(\infty)=(0,0)$ and maps $\infty$ to
$(w_0,w)$. It thus sends the geodesic line from $\infty$ to $(0,0)$
defined by $s\mapsto (0,0,e^{-2s+\ln 2})$ in horospherical
coordinates, hence by $s\mapsto [e^{-2s}:0:1]$ in homogeneous
coordinates by Equation \eqref{eq:horosphecoord}, to a geodesic line
from $(w_0,w)$ to $(0,0)$. An easy computation gives that this
geodesic line is
$$
s\mapsto \big(w_0(1+e^{2s}w_0)^{-1},\;w(1+e^{2s}w_0)^{-1}\big)\;.
$$
as wanted, after a time translation. \cqfd

\blemm\label{lem:disthoroghoro} For all $g\in U_q$ and $s>0$ such that
the horoballs $\H_s$ and $g\H_s$ have disjoint interiors, if $c_g$ is
the $(3,1)$-entry of the matrix $g$, then
$$
d(\H_s,g\H_s)=\frac 12\log\n(c_g)+\log\frac s2\,.
$$ 
\elemm 

\dem We follow \cite[Lem.~6.3]{ParPau10GT}. As seen in Section
\ref{sec:quathypspace}, if we had $c_g=0$, then $g$ would fix $\infty$
and would stabilise $\H_s$, which contradicts the assumption. Thus
$c_g\ne 0$.  Multiplying $g$ on the left and right by elements of
$\HHeis_{4n-1}$ does not change $c_g$ or $d(\H_s,g\H_s)$. We may
hence assume that $g(\infty)=(0,0)$ and $g^{-1}(\infty)=(0,0)$ (in the
coordinates $(w_0,w)$). Writing $g=\begin{pmatrix} a & \gamma^* &
b\\ \alpha & A & \beta \\ c & \delta^* & d\end{pmatrix}$, the first
condition implies that $a=0$ and $\alpha=0$, and the second one that
$\beta=0$ and $d=0$. The first and second equations of Formula
\eqref{eq:equationsUq} then imply that $\gamma=\delta=0$, the third
one implies that $A$ is unitary, and the fourth one gives $c\overline
b=1$. Thus,
$$
g=\begin{pmatrix} 0&0&\overline c^{-1}\\0&A&0\\c&0&0\end{pmatrix}
$$ 
with $A\in\Sp(n-1)$. It is easy to check, using the properties of
$\tr$ and $\n$, that
$$
g\H_s=\{(w_0,w)\in\HH\times\HH^{n-1}:\tr w_0-n(w)\ge s\n(c)\n(w_0)\}\,.
$$ 
The points of intersection of the geodesic line from $(0,0)$ to
$\infty$ and the horospheres $\partial\H_s$ (centred at $\infty$) and
$g\,\partial \H_s$ (centred at $(0,0)$) are $(\frac s2,0)$ and $(\frac
2{s\n(c)},0)$. The distance between them is as required by the 
statement.
\cqfd

\section{Measure computations in quaternionic hyperbolic spaces}
\label{sect:measurecomput}

Let $\Gamma$ be a nonelementary discrete group of isometries of
$\hnh$, let $\Lambda\Gamma$ be its limit set and let $\delta_\Ga$ be
its critical exponent.  In this section, we give proportionality
constants relating, on the one hand, Patterson, Bowen-Margulis and
skinning measures associated to some convex subsets and, on the other
hand, the corresponding Riemannian measures, in the quaternionic
hyperbolic case. These results were announced in
\cite[Chap.~7]{BroParPau19}, and we refer to Chapter 1 of op.~cit.~for
the background definitions and informations on the notions of this
section.

\medskip
We start by briefly recalling the construction of these measures. Let
$(\mu_{x})_{x\in \hnh}$ be a {\em Patterson density} for $\Ga$, that
is a family $(\mu_{x})_{x\in \hnh}$ of nonzero finite (nonnegative
Borel) measures on $\partial_{\infty}\hnh$ whose support is
$\Lambda\Ga$, such that $\ga_*\mu_x=\mu_{\ga x}$ and
$$
\frac{d\mu_{x}}{d\mu_{y}}(\xi)=e^{-\delta_{\Ga}\beta_{\xi}(x,\,y)}
$$
for all $\ga\in\Ga$, $x,y\in\hnh$ and (almost all)
$\xi\in\partial_{\infty}\hnh$.

For every $v\in T^1\hnh$, let $\pi(v)\in \hnh$ be its footpoint, and
let $v_-, v_+$ be the points at infinity of the geodesic line defined
by $v$. Let $x_0\in\hnh$ be a basepoint.  The {\em Bowen-Margulis
  measure} $\wt m_{\rm BM}$ for $\Ga$ on $T^1\hnh$ is defined, using
Hopf's parametrisation $v\mapsto (v_-,v_+,\beta_{v_+}(x_0,\pi(v))\,)$
from $T^1\hnh$ into $\partial_\infty \hnh\times\partial_\infty
\hnh\times\RR$, by
\begin{equation}\label{eq:defBM}
d\wt m_{\rm BM}(v)=e^{-\delta_{\Ga}(\beta_{v_{-}}(\pi(v),\,x_{0})+
\beta_{v_{+}}(\pi(v),\,x_{0}))}\;
d\mu_{x_{0}}(v_{-})\,d\mu_{x_{0}}(v_{+})\,dt\,.  
\end{equation}
Note that in the right hand side of this equation, $\pi(v)$ may be
replaced by any point $x'$ on the geodesic line defined by $v$, since
$\beta_{v_{-}}(\pi(v),\,x')+ \beta_{v_{+}}(\pi(v),\,x')=0$.  We will
use this elementary observation in the proof of Lemma
\ref{lem:computheisen} {\em (ii)}.  The measure $\wt m_{\rm BM}$ is nonzero,
independent of $x_{0}$, is invariant under the geodesic flow,
the antipodal map $v\mapsto -v$ and the action of $\Ga$. Thus, it
defines a nonzero measure $m_{\rm BM}$ on $T^1\Ga\bs\hnh$ which is
invariant under the geodesic flow of $\Ga\bs\hnh$ and the antipodal
map, called the {\em Bowen-Margulis measure} on $\Ga\bs \hnh$.

\medskip 
Let $D$ be a nonempty proper closed convex subset of $\hnh$, with
stabiliser $\Ga_{D}$ in $\Ga$, such that the family $(\ga
D)_{\ga\in\Ga/\Ga_{D}}$ is locally finite in $\hnh$. We denote by
$\normalpm D$ the {\it outer/inner unit normal bundle} of $\partial
D$, that is, the set of $v\in T^1\hnh$ such that $\pi(v)\in \partial
D$ and the closest point projection on $D$ of $v_\pm\in\partial_\infty
\hnh- \partial_\infty D$ is $\pi(v)$. Using the endpoint homeomorphism
$v\mapsto v_\pm$ from $\normalpm {D}$ to $\partial_{\infty}\hnh
-\partial_{\infty}D$, we defined in \cite{ParPau14ETDS} (generalising
the definition of Oh and Shah \cite[\S 1.2]{OhSha12} when $D$ is a
horoball or a totally geodesic subspace in $\hnr$) the outer/inner
{\em skinning measure} $\wt\sigma^\pm_{D}$ of $\Ga$ on $\normalpm{D}$,
by
\begin{equation}\label{eq:defskin}
d\wt\sigma^\pm_{D}(v) =  
e^{-\delta_\Ga\,\beta_{v_{\pm}}(\pi(v),\,x_{0})}\,d\mu_{x_{0}}(v_{\pm})\,.
\end{equation}
The measure $\wt\sigma^\pm_{D}$ is independent of $x_{0}$. It is
nonzero if $\Lambda\Ga$ is not contained in $\partial_{\infty}D$, and
it satisfies $\wt\sigma^\pm_{\ga D} =\ga_{*}\wt\sigma^\pm_{D}$ for
every $\ga\in\Ga$. The measure $\sum_{\ga \in \Ga/\Ga_{D}} \;\ga_*
\wt\sigma^\pm_{D}$ is a well defined \mbox{$\Ga$-invariant} locally
finite measure on $T^1\hnh$.  Hence, it induces a locally finite
measure $\sigma^\pm_{D}$ on $\Ga\bs T^1\hnh$, called the outer/inner
{\em skinning measure} of $D$ in $\Ga\bs T^1\hnh$. Note that if
$\iota:v\mapsto -v$ is the antipodal map, then $\iota_*\wt\sigma^-_{D}
=\wt\sigma^+_{D}$. In particular $\pi_*\wt\sigma^-_{D}
=\pi_*\wt\sigma^+_{D}$, and the measures $\sigma^-_{D}$ and
$\sigma^+_{D}$ have the same total mass.

We will denote the standard Lebesgue measures on the Euclidean spaces
$\HH^{n-1}$ and $\Im\;\HH$ by $d\zeta$ and $du$ respectively, so that
the usual left Haar measure $d\lambda_{4n-1}$ on the Lie group
$\HHeis_{4n-1}$ is
\begin{equation}\label{eq:defhaarlambda}
d\lambda_{4n-1}(\zeta,u) = d\zeta du\,.
\end{equation}
In horospherical coordinates, the volume form of $\hnh=
\HHeis_{4n-1} \times \;\mathopen{]}0,+\infty[$ is\footnote{See also
\cite[page 301]{KimPar03} with a different normalisation.}
\begin{equation}\label{eq:volhnh}
d\operatorname{vol}_{\hnh}(\zeta,u,t)=
\frac{1}{16\,t^{2n+2}}\;d\zeta\,du\,dt\;.
\end{equation}

We begin by giving a lemma that relates the Riemannian volume of a
Margulis cusp neighbourhood with the Riemannian volume of its
boundary, close to \cite[Lem.~3.1]{KimPar03}.

\blemm\label{lem:horobvol}
Let $D$ be a horoball in $\hnh$ and let $\Ga$ be a discrete group of
isometries of $\hnh$ preserving $\partial D$ (hence $D$). Then
$\Vol(\Ga\backslash\partial D)=(4n+2)\Vol(\Ga\backslash D)$.  
\elemm

\dem 
Since the group of isometries of $\hnh$ acts transitively on the set
of horospheres of $\hnh$, we may assume that $D=\H_1$.  The horosphere
centred at $\infty$ through a point $(\zeta,u,t) \in \H_1$ is equal to
$\partial \H_t$ and its orthogonal geodesic line at this point is
$s\mapsto (\zeta,u,e^{2s})$, hence
$$
d\vol_{\hnh}(\zeta,u,t)=d\vol_{\partial \H_t} (\zeta,u,t)
\,\frac{dt}{2t}\,.
$$ 
By Equation \eqref{eq:volhnh}, we hence have 
\begin{equation}\label{eq:volhorosphzetau}
d\vol_{\partial \H_t}(\zeta,u,t)=\frac{1}{8\,t^{2n+1}}\;d\zeta\,du
\end{equation} 
for every $t>0$, therefore $d\vol_{\partial \H_t}
(\zeta,u,t)=\frac{1}{t^{2n+1}}\; d\vol_{\partial \H_1} (\zeta,u,1)$. The
homeomorphism from $\partial \H_t$ to $\partial \H_1$ defined by
$(\zeta,u,t)\mapsto (\zeta,u,1)$ commutes with the action of $\Ga$.
Thus,
\begin{align*}
\Vol(\Ga\backslash\H_1) &=\int_{\Ga\backslash\H_1}
\dvol_{\hnh}(\zeta,u,t)=\int_{t=1}^{+\infty}\int_{\Ga\backslash\partial \H_t}
\dvol_{\partial \H_t} (\zeta,u,t)\,\frac{dt}{2t}
\\ & =\int_{t=1}^{+\infty}\int_{\Ga\backslash\partial \H_1}
\dvol_{\partial \H_1} (\zeta,u,1)\,\frac{dt}{2t^{2n+2}}=
\frac{1}{4n+2}\,\Vol(\Ga\backslash \partial \H_1)
\;. \;\;\;\Box
\end{align*}

\medskip
Let $\Ga$ be a lattice in $\Isom(\hnh)$, that is, a discrete group of
isometries of $\hnh$ such that the orbifold $\Ga\bs\hnh$ has finite
volume. Its critical exponent is
\begin{equation}\label{eq:critexpolat}
\delta_\Ga=4n+2
\end{equation} 
(see for instance \cite[Theo.~4.4 (i)]{Corlette90}). The Patterson
density $(\mu_x)_{x\in\hnh}$ of $\Ga$ is uniquely defined up to a
multiplicative constant, and is independent of $\Ga$. We will choose
the normalisation as follows. Let $\mu_\infty$ be the
$\HHeis_{4n-1}$-invariant measure on $\partial_\infty \hnh-\{\infty\}$
defined (see for instance \cite[Prop.~7.2]{BroParPau19}) by
\begin{equation}\label{eq:defmuinfty}
\mu_\infty=\lim_{t\ra+\infty}\; e^{\delta_{\Ga} t}\,\mu_{\rho(t)}
\end{equation} 
where $\rho$ is the geodesic ray starting from any point in $\partial
\H_1$ and converging to $\infty$. By the uniqueness property of Haar
measures on $\HHeis_{4n-1}$, we may uniquely normalise the Patterson
density so that $\mu_\infty$ coincides with $\lambda_{4n-1}$ on
$\partial_\infty \hnh-\{\infty\}= \HHeis_{4n-1}$, that is
$$
d\mu_\infty(\xi,r)= d\lambda_{4n-1}(\xi,r) =d\xi\,dr\;.
$$

The various computations of Patterson, Bowen-Margulis and skinning
measures are gathered in the following statement.

\blemm \label{lem:computheisen} Let $\Ga$ be a lattice in
$\Isom(\hnh)$, and let $(\mu_x)_{x\in\hnh}$ be its Patterson density,
normalised as above.  For all $x=(\zeta,u,t)$ and $x'=(\zeta',u',t')$
in $\hnh$, for all $(\xi,r)$ in $\partial_\infty \hnh-\{\infty\}$ and
$v$ in $T^1\hnh$ such that $v_\pm\neq\infty$, we have

\smallskip\noindent
(i)~ 
$\displaystyle d\mu_x(\xi,r)=
\frac{t^{2n+1}}{d_{\rm Cyg}(x,(\xi,r))^{8n+4}}\;d\xi\,dr\;;$

\medskip\noindent
(ii) using a Hopf parametrisation $v\mapsto (v_-,v_+,s)$, 
$$
d\wt m_{\rm BM}(v)=
\frac{d\lambda_{4n-1}(v_-)\,d\lambda_{4n-1}(v_+)\,ds} 
{d_{\rm  Cyg}(v_-,v_+)^{8n+4}}\,;
$$

\noindent(iii) 
$$
\wt m_{\rm BM}= \frac{1}{2^{4n-4}}\; \vol_{T^1\hnh}\,,
$$ 
and in particular, if $M=\Ga\bs\hnh$, the total mass of the
Bowen-Margulis measure of $\Ga\bs T^1 \hnh$ is
$$
\|m_{\rm BM}\|= \frac{\pi^{2n}}{2^{4n-5}\,(2n-1)!}\;\Vol(M)\;;
$$

\medskip\noindent(iv) using the homeomorphism $v\mapsto v_+$ from
$\normalout\H_1$ to $\partial_\infty \hnh-\{\infty\}= \HHeis_{4n-1}$,
we have
$$
d\wt\sigma^+_{\H_1}(v)=d\lambda_{4n+1}(v_+)\,;
$$
for every horoball $D$ in $\hnh$, we have
$$
\pi_*\wt\sigma^\pm_{D}=8\,\vol_{\partial D}\,,
$$ 
and the total mass of the skinning measure of $D$ in $\Ga\bs T^1
\hnh$ is
$$
\|\sigma^\pm_{D}\|= 16(2n+1)\, \Vol(\Ga_{D}\bs D)\;;
$$

\medskip\noindent 
(v) for every geodesic line $D$ in $\hnh$, we have
$$
d\pi_*\wt \sigma^\pm_{D}=\frac{2n+1}{2^{4n-2}\,(4n-1)}\;
d\pi_*\vol_{\normalpm D}
$$
and, with $m$ the order of the pointwise stabiliser of $D$ in $\Ga$,
$$
\|\sigma^\pm_{D}\|= 
\frac{\pi^{2n-1}\,(2n+1)!}{m\,n\,(4n-1)!}\;
\Vol(\Ga_{D}\bs D)\,;
$$

\medskip\noindent 
(vi) for every quaternionic geodesic line $D$ in $\hnh$, we have
$$
d\pi_*\wt \sigma^\pm_{D}= \frac{1}{2^{4n-1}}\;d\pi_*\vol_{\normalpm D}
$$
and, with $m$ the order of the pointwise stabiliser of $D$ in $\Ga$,
$$
\|\sigma^\pm_{D}\|= 
\frac{\pi^{2n-2}}{m\,2^{4n-2}\,(2n-3)!}\;\Vol(\Ga_{D}\bs D)\,.
$$  
\elemm

\dem In the computations below, it is useful to note that Lemma
\ref{lem:busemann} implies that
\begin{equation}\label{eq:rewritelem7}
e^{-(4n+2)\,\beta_{(\xi,\,r)}(x,\,x')}= 
\frac{t^{2n+1}\;d_{\rm Cyg}(x',(\xi,r))^{8n+4}}
{(t')^{2n+1}\;d_{\rm Cyg}(x,(\xi,r))^{8n+4}}\,.
\end{equation}

\medskip\noindent {\em (i)~} 
The geodesic line from $(\xi,r)$ to $\infty$ goes through $\partial
\H_1$ at the point $(\xi,r,1)$. For all $\eta\in\partial_\infty
\hnh-\{\infty\}$, let $x_{\H_1,\,\eta}$ be the intersection point with
$\partial \H_1$ of the geodesic line from $\eta$ to $\infty$. By the
normalisation of $d\mu_\infty$ and by the definition of $\mu_\infty$
and the Radon-Nikodim property of the Patterson density, we have
$$ 
\frac{d\mu_x}{d\mu_\infty}(\eta)=
e^{-\delta_{\Ga}\beta_{\eta}(x,\,x_{\H_1,\,\eta})}\;,
$$ 
for all $x\in\hnh$ and (almost all) $\eta\in\partial_\infty
\hnh-\{\infty\}$. Hence we have
$$
\frac{d\mu_x}{d\xi\,dr} (\xi,r)=
\frac{d\mu_x}{d\mu_\infty}(\xi,r)
=e^{-\delta_{\Ga}\,\beta_{(\xi,\,r)}(x,\,(\xi,\,r,\,1))}\,.
$$ 
The result then follows from Equations \eqref{eq:critexpolat},
\eqref{eq:rewritelem7} and \eqref{eq:defidistCyg}.

\medskip\noindent {\em (ii)~} Note that if $x'$ is on the geodesic
line $\ell$ defined by $v$, then by Lemma \ref{lem:calcgeodline} and
the fact that $\ell$ is asymptotic near $v_-$ to the geodesic line
from $v_-$ to $\infty$, an easy computation using Equations
\eqref{eq:defidistCyg} and \eqref{eq:horosphecoord} shows that $d_{\rm
  Cyg} (x',v_-)^2\sim t'$ as $x'\ra v_-$.  Hence, by Equation
\eqref{eq:defBM} and the comment following it, by Equations
\eqref{eq:critexpolat} and \eqref{eq:rewritelem7}, by Assertion {\em
  (i)}, and by letting $x'$ converge to $v_-$ on the geodesic line
defined by $v$, we have

\begin{align*}
&d\wt m_{\rm BM}(v) 
=e^{-(4n+2)(\beta_{v_{-}}(x',\,x)+ \beta_{v_{+}}(x',\,x))}\; 
d\mu_{x}(v_{-})\,d\mu_{x}(v_{+})\,ds 
\\ &
=\Big(\frac{t'\,d_{\rm Cyg}(x,v_-)^{4}\,t'\,d_{\rm Cyg}(x,v_+)^{4}\,t^2}
{t\,d_{\rm Cyg}(x',v_-)^{4}\,t\,d_{\rm Cyg}(x',v_+)^{4}\;
d_{\rm Cyg}(x,v_-)^{4}\,d_{\rm Cyg}(x,v_+)^{4}}\Big)^{2n+1}\;
\\ &\;\;\;\;\;\;d\lambda_{4n-1}(v_{-})\,d\lambda_{4n-1}(v_{+})\,ds
\\ &
=\frac{1}{d_{\rm Cyg}(v_-,v_+)^{8n+4}}\;
d\lambda_{4n-1}(v_{-})\,d\lambda_{4n-1}(v_{+})\,ds\,.
\end{align*}

\smallskip\noindent 
{\em (iii)~} Recall that the Liouville measure $\vol_{T^1\hnh}$ (which
is the Riemannian measure for Sasaki's metric on $T^1\hnh$)
disintegrates under the fibration $\pi:T^1\hnh \ra \hnh$ over the
Riemannian measure $\vol_{\hnh}$ of $\hnh$, with conditional measures
the spherical measures on the unit tangent spheres:
$$
d\vol_{T^1\hnh}(v)=\int_{x\in \hnh} \dvol_{T^1_x\hnh}(v)\dvol_{\hnh}(x)\;.
$$ 
Let $x=(\zeta,u,t)\in \hnh$. Since the group $I_x$ of isometries of
$\hnh$ fixing $x$ acts transitively on $T^1_x\hnh$, since both $\mu_x$
and the Riemannian measure $\vol_{T^1_x\hnh}$ are invariant under
$I_x$, using the $I_x$-equivariant homeomorphism $v\mapsto v_+$ from
$T^1_x\hnh$ to $\partial_\infty\hnh$, we have, for all $v\in
T^1_x\hnh$ such that $v_+\neq \infty$, using Assertion {\em (i)} for
the last equality,
\begin{equation}\label{eq:relatvolfibunilambda}
d\vol_{T^1_x\hnh}(v)=\frac{\Vol(\SSS^{4n-1})}{\|\mu_x\|}\;d\mu_x(v_+)
=\frac{\Vol(\SSS^{4n-1})\;t^{2n+1}}{\|\mu_x\|\,d_{\rm Cyg}(x,v_+)^{8n+4}}
\;d\lambda_{4n-1}(v_+)\,.
\end{equation}
By homogeneity, by Assertion {\em (i)}, by Equation
\eqref{eq:defidistCyg} with $ (\zeta,u,t)=(0,0,1)$ and $(\zeta',u',t')
=(\xi,r,0)$, by using the spherical coordinates in the Euclidean
spaces $\HH^{n-1}$ and $\Im\;\HH$ of real dimensions $4n-4$ and $3$ so
that $d\xi= s^{4n-5}ds\,d\vol_{\SSS^{4n-5}}$ and $dr=\rho^2 d\rho\,
d\vol_{\SSS^2}$, and by using the changes of variables $\rho\mapsto
\frac{\rho}{s^2+1}$ and $s\mapsto s^2$, we have
\begin{align*}
\|\mu_x\|&=\|\mu_{(0,0,1)}\|= \int_{\HH^{n-1}\times
  \Im\;\HH}\frac{d\xi\,dr}{((\n(\xi)+1)^2+\n(r))^{2n+1}} \\ & =
\Vol(\SSS^{4n-5})\Vol(\SSS^2) \iint_{0}^{+\infty}
\frac{s^{4n-5}\rho^{2}\;ds\,d\rho}{((s^2+1)^2+\rho^2)^{2n+1}}
\\ & = \pi\,\Vol(\SSS^{4n-5})
\int_{-\infty}^{+\infty}\frac{\rho^{2}\;d\rho}{(1+\rho^2)^{2n+1}}\;
\int_0^{+\infty}\frac{\;s^{2n-3}\,ds}{(s+1)^{4n-1}}\;.
\end{align*}
By the residue formula at a pole of order $2n+1$ and by Leibniz
formula, considering the map $f:z\mapsto \frac{1}{(z+i)^{2n+1}}$, we
have
\begin{align*}
&\int_{-\infty}^{+\infty}\frac{\rho^2\;d\rho}{(\rho^2+1)^{2n+1}}=
2i\pi\Res_{z=i} \frac{z^2}{(z^2+1)^{2n+1}}=
2i\pi\frac{1}{(2n)!}\frac{\partial^{2n}}{\partial z^{2n}}\Big|_{z=i} 
(z^2f(z))\\ =&
\;\frac{2i\pi}{(2n)!}\Big(z^2\,\frac{\partial^{2n}f}
{\partial z^{2n}}+4n\,z\,\frac{\partial^{2n-1}f}
{\partial z^{2n-1}}+2n(2n-1)\frac{\partial^{2n-2}f}
{\partial z^{2n-2}}\Big)_{z=i}=
\frac{\pi\,n\,(4n-2)!}{2^{4n-2}((2n)!)^2}\;.
\end{align*}
By integration by part and by induction, we have
$$
\int_0^{+\infty}\frac{s^{2n-3}\;ds}{(s+1)^{4n-1}}=
\frac{(2n-3)!}{(4n-2)\dots (2n+2)}
\int_0^{+\infty}\frac{ds}{(s+1)^{2n+2}}
=\frac{(2n-3)!(2n)!}{(4n-2)!}\;.
$$
Since $\Vol(\SSS^{4n-1})= \frac{\pi^2}{(2n-1)(2n-2)}
\,\Vol(\SSS^{4n-5})$, we hence have
$$
\|\mu_x\|=\frac{1}{2^{4n-1}}\;\Vol(\SSS^{4n-1})\;.
$$
Hence, by Equations \eqref{eq:volhnh} and
\eqref{eq:relatvolfibunilambda}, using the homeomorphism $v\mapsto
(v_+,\pi(v)=(\zeta,u,t))$ from $T^1\hnh$ to $\partial_\infty \hnh
\times \hnh$, we have, for all $v\in T^1\hnh$ such that $v_+\neq
\infty$,
\begin{equation}\label{eq:liouville}
d\vol_{T^1\hnh}(v) =
\frac{2^{4n-5}}{t\,d_{\rm Cyg}((\zeta,u,t),v_+)^{8n+4}}
\;d\lambda_{4n-1}(v_+)\,d\zeta\,du\,dt\,.
\end{equation}

\medskip
Now, let us consider the map $F:\HHeis_{4n-1}\times\RR\,\ra\,\hnh$
defined by
\begin{align*}
(\xi,r,s)\mapsto  \Big(\;&\zeta=\xi\;(1+(\n(\xi)+r)\,e^{2s})^{-1},\; \\ &
u=\Im\;\big((\n(\xi)+r)(1+(\n(\xi)+r)\,e^{2s})^{-1}\big),\;\\ &
t=\frac{(\n(\xi)^2+\n(r))\,e^{2s}}{\n(1+(\n(\xi)+r)\,e^{2s})}\;\Big)\;.
\end{align*}
Note that $F(0,i,0)=(0,\frac{i}{2},\frac{1}{2})$. By Lemma
\ref{lem:calcgeodline} and Equation \eqref{eq:horosphecoord}, the map
$s\mapsto F(\xi,r,s)$ is a geodesic line in $\hnh$ starting from
$(\xi,r)$ and ending at $(0,0)$. On this geodesic, $s$ and the time
parameter in Hopf's parametrisation differ only by an additive
constant, hence have the same differential.

Recall that by homogeneity, the two measures $\wt m_{\rm BM}$ and
$\vol_{T^1\hnh}$ are proportional. Hence, computing their (constant)
Radon-Nikodym derivative at $v\in T^1\hnh$ such that $v_-=(0,i)$ and
$\pi(v)=(0,\frac{i}{2},\frac{1}{2})$ (so that $v$ is tangent to the
geodesic line $s\mapsto F(0,i,s)$ at $s=0$, hence $v_+=(0,0)$), we
have, by Assertion {\em (ii)~} with $(\xi,r)$ parametrising $v_-$ and
by Equations  \eqref{eq:liouville} and \eqref{eq:defidistCyg},
\begin{align*}
\frac{d\vol_{T^1\hnh}}{d\wt m_{\rm BM}}&= 
\frac{2^{4n-5}\,d_{\rm Cyg}((0,i,0),(0,0,0))^{8n+4}} 
{\frac 12\,d_{\rm Cyg}((0,\frac{i}{2},\frac{1}{2}),(0,0,0))^{8n+4}}\;
\frac{d\zeta\,du\,dt}{d\xi\,dr\,ds}(0,i,0)\\ &
=2^{6n-3}\;
\frac{d\zeta\,du\,dt}{d\xi\,dr\,ds}(0,i,0)
\,.
\end{align*}
Let us compute the Jacobian at $(0,i,0)$ of the map $F:(\xi,r,s)
\mapsto (\zeta,u,t)$. At the point $(0,i,0)$, we have, using the
canonical basis $i,j,k$ of $\Im\;\HH$ in order to write
$r=r_1i+r_2j+r_3k$ and $u=u_1i+u_2j+u_3k$,
$$
\frac{\partial \zeta}{\partial \xi}=
\frac{1}{1+i}\operatorname{Id}_{\HH^{n-1}},\;\;\;
\frac{\partial \zeta}{\partial r}=
\frac{\partial \zeta}{\partial s}=0
,\;\;\;\frac{\partial t}{\partial s}=
\frac{\partial t}{\partial \xi} =0\;,\;\;\;
\frac{\partial u}{\partial \xi}= 0\;,
$$
$$
\frac{\partial u}{\partial r}=\begin{pmatrix}0 & 0 & 0\\
0 & \frac{1}{2} & 0\\ 0 & 0 &\frac{1}{2} \end{pmatrix},\;\;\;
\frac{\partial u}{\partial s}=
\begin{pmatrix}-1 \\ \;0 \\\; 0 \end{pmatrix},\;\;\;
\frac{\partial t}{\partial r}=
\begin{pmatrix}\frac{1}{2} & 0 & 0 \end{pmatrix}\;.
$$ 
Since the Jacobian matrix of $F$ at $(0,i,0)$ is block diagonal when
the variables are separated into the $4(n-1)$ first ones and the last
$4$ ones, since the multiplication by $\frac{1}{1+i}$ in $\HH^{n-1}$
is a Euclidean homothety of ratio $\frac{1}{\sqrt{2}}$, and since the
determinant of the $4\times 4$ matrix of the partial derivatives of
$u,t$ with respect to $r,s$ has absolute value $\frac18$, the Jacobian
of $F$ at $(0,i,0)$ is equal to $\big(\frac{1}{\sqrt{2}}\big)^{4(n-1)}
\frac18 =\frac{1}{2^{2n+1}}$.  The first claim of Assertion {\em
  (iii)} follows

The second claim follows from the facts that $\Vol(T^1M)=
\Vol(\SSS^{4n-1})\,\Vol (M)$ and that $\Vol(\SSS^{4n-1})=
\frac{2\,\pi^{2n}}{(2n-1)!}$.

\medskip\noindent {\em (iv)~} By the definition of the skinning
measure $\wt\sigma^+_{\H_1}$ in Equation \eqref{eq:defskin} and of the
measure $\mu_\infty$ in Equation \eqref{eq:defmuinfty}, we have
$$
d\wt\sigma^+_{\H_1}(v) = d\mu_{\infty}(v_{+})
$$ 
for every $v\in \normalout \H_1$, since $\beta_{v_{+}}(\pi(v),\rho(t))
=-t+\smallo(1)$ as $t\ra+\infty$. The first claim of Assertion {\em
(iv)} follows by the normalisation of the Patterson density.

By Equation \eqref{eq:volhorosphzetau}, we have
\begin{equation}\label{eq:volHun}
d\vol_{\partial \H_1}(\zeta,u,1)= \frac{1}{8}\;d\zeta\,du\,.
\end{equation}
Hence $\pi_*\wt \sigma^\pm_{\H_1}= 8\, \vol_{\partial \H_1}$, and by the
transitivity of the isometry group of $\hnh$ on the set of horoballs
in $\hnh$, the second claim of Assertion {\em (iv)} follows.
Therefore, by Lemma \ref{lem:horobvol},
$$
\|\sigma^\pm_{D}\|=\|\pi_*\sigma^\pm_{D}\|=
8\,\Vol(\Ga_{D}\bs\partial D)= 16(2n+1)\,\Vol(\Ga_{D}\bs D)\,.
$$

\medskip\noindent 
{\em (v)~} By the transitivity of the isometry group of $\hnh$ on the
set of its geodesic lines, we may assume that $D$ is the geodesic
line in $\hnh$ with points at infinity $(0,0)$ and $\infty$.  The map
from the full-measure open subset $\big\{(\zeta,u)
\in\HHeis_{4n-1}\,:\;\zeta\neq0,\;u\neq0\big\}$ in $\HHeis_{4n-1}$ to
the product manifold $\SSS^{4n-5}\times\SSS^{2} \times\,]0,+\infty[
    \,\times\, ]0, \frac{\pi}{2}[$ defined by
\begin{equation}\label{eq:coordsigmarhotheta}
(\zeta,u)\mapsto \Big(\sigma=\frac{\zeta}{\n(\zeta)^{\frac12}},\;
w=\frac{u}{\n(u)^{\frac12}},\;
\rho=(\n(\zeta)^2+\n(u))^{1/2},\;
\theta=\arctan\frac{\n(u)^{1/2}}{\n(\zeta)}\Big)
\end{equation}
is a diffeomorphism. Since 
\begin{equation}\label{eq:inversepolar}
\n(\zeta)=\rho\cos\theta\;\;\; {\rm and}
\;\;\;\n(u)^{1/2}=\rho\sin\theta\,,
\end{equation} we have
\begin{align}
d\zeta\,du&=\frac{1}{2} 
\n(\zeta)^{2n-3}\,d(\n(\zeta))\dvol_{\SSS^{4n-5}}\Big(\frac{\zeta}
{\n(\zeta)^{1/2}}\Big)\,\n(u)\,d(\n(u)^{1/2})\,\dvol_{\SSS^{2}}
\Big(\frac{u}{\n(u)^{1/2}}\Big)\nonumber
\\ & =\frac{1}{2} \cos^{2n-3}\theta\;\sin^2\theta\;\rho^{2n}\,
\dvol_{\SSS^{4n-5}}(\sigma)\,\dvol_{\SSS^{2}}(w)\; d\rho\;d\theta\,.
\label{eq:measinpolar}
\end{align}
Using respectively in the following sequence of equalities

\smallskip$\bullet$~ the definition of the skinning measure in
Equation \eqref{eq:defskin} with basepoint $x_0=(0,0,1)$ and the
homeomorphism sending $v\in\normalout D$ to $v_+=(\zeta,u)\in
\HHeis_{4n-1} -\{(0,0)\}$, Equation \eqref{eq:critexpolat} and Lemma
\ref{lem:orthprojgeod},

\smallskip$\bullet$~ Equation \eqref{eq:rewritelem7} and Assertion
          {\em (i)},

\smallskip$\bullet$~ 
Equation \eqref{eq:defidistCyg}, and 

\smallskip$\bullet$~ Equations \eqref{eq:inversepolar} and
\eqref{eq:measinpolar},

\smallskip\noindent 
we have
\begin{align*}
d\wt \sigma^+_{D}(v)&= 
e^{-(4n+2)\,\beta_{(\zeta,u)}(\,(0,\,0,\,(\n(\zeta)^2+\n(u))^{1/2}),\;(0,\,0,\,1)\,)}\;
d\mu_{(0,\,0,\,1)}(\zeta,u)\\ & =\frac{(\n(\zeta)^2+\n(u))^{(2n+1)/2}}
{d_{\rm Cyg}((0,\,0,\,(\n(\zeta)^2+\n(u))^{1/2}),(\zeta,u,0))^{8n+4}}
\;d\zeta \,du \\ & = \Big(\frac{(\n(\zeta)^2+\n(u))^{1/2}}
{(\n(\zeta)+(\n(\zeta)^2+\n(u))^{1/2})^2+\n(u)}\Big)^{2n+1}\;
d\zeta \,du \\ &=
\frac{\cos^{2n-3}\theta\;\sin^2\theta}{2^{2n+2}\,(1+\cos\theta)^{2n+1}}
\dvol_{\SSS^{4n-5}}(\sigma)\,\dvol_{\SSS^{2}}(w)\;
\frac{d\rho}{\rho}\;d\theta\,.
\end{align*}
Thus,
\begin{equation}\label{eq:formskinningmeasgeod}
d\pi_*\wt \sigma^+_{D}(0,0,\rho)=
\frac{c'_n\,\Vol(\SSS^{4n-5})\Vol(\SSS^{2})}{2^{2n+2}}\;
\frac{d\rho}{\rho}\;,
\end{equation}
where, using the change of variable $t=\tan\frac{\theta}{2}$,
$$
c'_n=\int_{0}^{\frac{\pi}{2}}\frac{\cos^{2n-3}\theta\;\sin^2\theta}
{(1+\cos\theta)^{2n+1}}\;d\theta
=\frac{1}{2^{2n-2}}\int_{0}^1(1-t^2)^{2n-3}\,t^2\,(1+t^2)\;dt\;.
$$
With $I_{p,q}=\int_{-1}^1 t^{2p}(1-t^2)^q\;dt$, we have by integration by
part and by induction
$$
I_{p,q}= 
\frac{2^{2q+1}\,q!\,(2p)!\;(p+q)!}{p!\;(2p+2q+1)!}\;.
$$
Hence $c'_n=\frac{1}{2^{2n-1}}\;(I_{1,\,2n-3}+I_{2,\,2n-3})=
{\displaystyle \frac{2^{2n-1}\,(2n-3)!\,(2n-1)!\,(2n+1)}{(4n-1)!}}$\;.

\medskip 
The next step is to obtain an expression similar to Equation
\eqref{eq:formskinningmeasgeod} for the Riemannian measure of the
submanifold $\normalout D$ of $T^1\hnh$ (endowed with Sasaki's
metric). For every $x\in D$, let us denote by $\nu^1_x D$ the fiber
over $x$ of the normal bundle map $v\mapsto \pi(v)$ from $\normalout
D$ to $D$. We endow $\nu^1_x D$ with the spherical metric induced by
the scalar product of the tangent space $T_x\hnh$ at $x$. The
Riemannian measure of $\normalout D$ disintegrates under this
fibration over the Riemannian measure of $D$ as
$$
d\vol_{\normalout D}(v)=\int_{x\in D}d\vol_{\nu^1_x D}(v)\dvol_{D}(x)\,.
$$
By looking at the expression \eqref{eq:riemethorosphecoord} of the
Riemannian metric of $\hnh$ in horospherical coordinates, using the
homeomorphism $\rho \mapsto x=(0,0,\rho)$ from $]0,+\infty[$ to $D$,
we have
$$
d\vol_{D}(x)=\frac{d\rho}{2\rho}\,.
$$
Hence
\begin{equation}\label{eq:formvolnormbundmeasgeod}
d\pi_*\vol_{\normalout D}(0,0,\rho)=
\Vol(\SSS^{4n-2})\;\frac{d\rho}{2\rho}\,.
\end{equation}
We have $\Vol(\SSS^{4n-2})
=\frac{2^{4n-1}\,\pi^{2n-1}\,(2n-1)!}{(4n-2)!}$ and
$\Vol(\SSS^{4n-5})=\frac{2\;\pi^{2n-2}} {(2n-3)!}$.  
Equations \eqref{eq:formskinningmeasgeod} and
\eqref{eq:formvolnormbundmeasgeod} give the first claim of Assertion
{\em (v)}.

The second one follows, since pushforwards of measures preserve their
total mass, and since $\Vol(\Ga_{\normalpm D}\bs \normalpm D)=
\frac{\Vol(\SSS^{4n-2})}{m}\;\Vol(\Ga_{D}\bs D)$.

\bigskip\noindent {\em (vi)~} By the transitivity of the isometry
group of $\hnh$ on the set of its quaternionic geodesic lines, we may
assume that $D$ is the quaternionic geodesic line
$C=\{(w_0,w)\in\hnh\;:\;w=0\}$ or, in horospherical coordinates,
$C=\{(\zeta,u,t)\in \hnh\;:\; \zeta =0\}$.

\medskip
Hence, using the homeomorphism from $\normalout C$ to $\{(\zeta,u)
\in\HHeis_{4n-1}\;:\;\zeta\neq 0\}$ sending a normal unit vector $v$ to
its point at infinity $v_+=(\zeta,u)$, by the definition of the
skinning measure in Equation \eqref{eq:defskin} with basepoint
$x_0=(0,0,1)$, by Equation \eqref{eq:critexpolat}, by Lemma
\ref{lem:projsurquatgeodline}, by Equations \eqref{eq:rewritelem7} and
\eqref{eq:defidistCyg}, and by Assertion {\em (i)}, we have
\begin{align*}
d\wt \sigma_{C}(v)&= 
e^{-(4n+2)\,\beta_{(\zeta,u)}(\,(0,\,u,\,\n(\zeta)),\;(0,\,0,\,1)\,)}\;
d\mu_{(0,\,0,\,1)}(\zeta,u)=\frac{1}{2^{4n+2}\n(\zeta)^{2n+1}}\;d\zeta\,du\\ &=
\frac{1}{2^{4n+3}\n(\zeta)^{4}}\;
d(\n(\zeta))\dvol_{\SSS^{4n-5}}\Big(\frac{\zeta}
{\n(\zeta)^{1/2}}\Big)\,du\,.
\end{align*}
In particular,
$$
d\pi_*\wt \sigma_{C}(0,u,\n(\zeta)) =
\frac{\Vol(\SSS^{4n-5})}{2^{4n+3}}\;du\;\frac{d(\n(\zeta))}{\n(\zeta)^{4}}\;.
$$

For every $x\in C$, let us denote by $\nu^1_x C$ the fiber over $x$ of
the normal bundle map $v\mapsto \pi(v)$ from $\normalout C$ to $C$,
endowed with the spherical metric induced by the scalar product of the
tangent space $T_x\hnh$ at $x$. The Riemannian measure of $\normalout
C$ disintegrates under this fibration over the Riemannian measure of
$C$ as
$$
d\vol_{\normalout C}(v)=\int_{x\in C}d\vol_{\nu^1_x C}(v)\dvol_{C}(x)\,.
$$ 
Using Equation \eqref{eq:riemethorosphecoord} and the homeomorphism
$(u,t=\n(\zeta)) \mapsto x=(0,u,t)$ from $\Im\;\HH
\times\,]0,+\infty[$ to $C$, we have
$$
d\vol_{C}(x)=\big(\frac{1}{2\,t}\big)^4du\,dt=
\frac{1}{2^4}\;du\;\frac{d(\n(\zeta))}{\n(\zeta)^4}\,.
$$
Hence
$$
d\pi_*\vol_{\normalout C} (x)=\Vol(\SSS^{4n-5})\dvol_{C}(x)=
\frac{\Vol(\SSS^{4n-5})}{2^4}\;du\;\frac{d(\n(\zeta))}{\n(\zeta)^4}\;.
$$
The result follows as in the end of the proof of the previous
Assertion. 
\cqfd

\section{Equidistribution and counting in quaternionic hyperbolic 
geometry}
\label{sect:mertens}

In this section, we first use the general results of
\cite{ParPau17ETDS} (see also \cite{BroParPau19}) and the computations
of Section \ref{sect:measurecomput} to give explicit asymptotic
counting and equidistribution results on the number of common
perpendiculars that are shorter than a given bound between two
properly embedded locally convex proper closed subsets of
$\Ga\bs\hnr$, for any lattice $\Ga$ in $\PU_q$. Using Sections
\ref{sec:cuspnumb} and \ref{sec:computvol}, we then give two
arithmetic applications, generalising Theorems \ref{theo:countintro}
and \ref{theo:equidisintro} in the introduction. We refer to
\cite{ParPau20b} for geometric applications.

\medskip
Let $\Ga$ be a lattice in $\PU_q$. Let $D^-$ and $D^+$ be nonempty
proper closed convex subsets of $\hnh$, with stabilisers $\Ga_{D^-}$
and $\Ga_{D^+}$ in $\Ga$ respectively, such that the families $(\ga
D^-)_{\ga\in\Ga/\Ga_{D^-}}$ and $(\ga D^+)_{\ga\in\Ga/\Ga_{D^+}}$ are
locally finite in $\hnh$. With the measures defined at the beginning
of Section \ref{sect:measurecomput}, let
$$
c(D^-,D^+)=\frac{\|\sigma^+_{D^-}\|\;\|\sigma^-_{D^+}\|}
{\delta_\Ga\;\|m_{\rm BM}\|}\;.
$$

For all $\ga,\ga'$ in $\Ga$, the convex sets $\ga D^-$ and $\ga' D^+$
have a common perpendicular if and only if their closures
$\overline{\ga D^-}$ and $\overline{\ga' D^+}$ in $\hnh\cup
\partial_\infty \hnh$ do not intersect.  We denote by
$\alpha_{\ga,\,\ga'}$ this common perpendicular, starting from $\ga
D^-$ at time $t=0$, and by $\ell(\alpha_{\ga,\,\ga'})$ its length.
The {\em multiplicity} of $\alpha_{\ga,\ga'}$ is
$$
m_{\ga,\ga'}=
\frac 1{\card(\ga\Ga_{D^-}\ga^{-1}\cap\ga'\Ga_{D^+}{\ga'}^{-1})}\,,
$$ 
which equals $1$ for all $\ga,\ga'\in\Ga$ when $\Ga$ acts freely on
$T^1\hnh$ (for instance when $\Ga$ is torsion-free). For all $s> 0$
and $x\in\partial D^-$, let
$$
m_s(x)=\sum_{\ga\in \Ga/\Ga_{D^+}\;:\;
  \overline{D^-}\,\cap \,\overline{\ga D^+}\,= \emptyset,\;
\alpha_{e,\, \ga}(0)=x,\; \ell(\alpha_{e,\, \ga})\leq s} m_{e,\ga}
$$ 
be the multiplicity of $x$ as the origin of common perpendiculars
with length at most $t$ from $D^-$ to the elements of the $\Ga$-orbit
of $D^+$. For every $s> 0$, let
$$
\N_{D^-,\,D^+}(s)=
\sum_{
(\ga,\,\ga')\in \Ga\bs((\Ga/\Ga_{D^-})\times (\Ga/\Ga_{D^+}))\;:\;
\overline{\ga D^-}\,\cap \,\overline{\ga' D^+}\,=\emptyset,\; 
\ell(\alpha_{\ga,\, \ga'})\leq s} m_{\ga,\ga'}
\;,
$$
where $\Ga$ acts diagonally on $\Ga\times\Ga$. When $\Ga$ has no
torsion, $\N_{D^-,\,D^+}(s)$ is the number (with multiplicities coming
from the fact that $\Ga_{D^\pm}\bs D^\pm$ is not assumed to be
embedded in $\Ga\bs\hnh$) of the common perpendiculars of length at
most $s$ between the images of $D^-$ and $D^+$ in $\Ga\bs\hnh$.

\medskip
Let us determine some constants before stating Theorem 8.1 giving the
asymptotics of $\N_{D^-,\,D^+}(s)$ as $s\ra +\infty$, and its
associated equidistribution claim. We assume from now on that $D^-$ is
a horoball in $\hnh$ centred at a parabolic fixed point of $\Ga$. We
assume from now on that $D^+$ is one of the following three
possibilities, we denote by $m^+$ the cardinality of the pointwise
stabiliser of $D^+$ in $\Ga$ and we compute $c(D^-,D^+)$ using Lemma
\ref{lem:computheisen}. If $D^+$ is also a horoball in $\hnh$ centred
at a parabolic fixed point of $\Ga$, then
$$
c(D^-,D^+)=\frac{2^{4n+1}\,(2n+1)!}{n\;\pi^{2n}}\;
\frac{\Vol(\Ga_{D^-}\bs D^-)\Vol(\Ga_{D^+}\bs D^+)}{\Vol(\Ga\bs\hnh)}\,.
$$
If $D^+$ is a geodesic line in $\hnh$ such that $\Ga_{D^+}\bs D^+$ is
compact, then
$$
c(D^-,D^+)=\frac{2^{4n}\,(2n-1)!(2n+1)!}{\pi\,m^+\,(4n)!}
\frac{\Vol(\Ga_{D^-}\bs D^-)\Vol(\Ga_{D^+}\bs D^+)}
{\Vol(\Ga\bs\hnh)}\,.
$$
If $D^+$ is a quaternionic geodesic line in $\hnh$ such that
$\Ga_{D^+}\bs D^+$ has finite volume, then
$$
c(D^-,D^+)=\frac{2\,(n-1)\,(2n-1)}{\pi^2\,m^+}
\frac{\Vol(\Ga_{D^-}\bs D^-)\Vol(\Ga_{D^+}\bs D^+)}
{\Vol(\Ga\bs\hnh)}\,.
$$
Lemma \ref{lem:computheisen} (iv) also gives that
$$
\frac{1}{\|\sigma^-_{D^-}\|}\;d\,\pi_*\sigma^+_{D^-}= 
\frac{1}{2\,(2n+1)\,\Vol(\Ga_{D^-}\bs D^-)}\;d\vol_{\partial D^-}\;.
$$

Recall that every lattice in $\PU_q$ is arithmetic, by the works of
Margulis, Corlette, Gromov-Schoen, see \cite[Theo.~8.4]{GroSch92}. The
following counting and equidistribution result of common
perpendiculars follows from \cite[Theo.~15 (2)]{ParPau17ETDS} (with
the remark preceding it concerning the proof by Kleinbock-Margulis and
Clozel of the exponential mixing property for the Sobolev regularity
of the geodesic flow), see also \cite[\S 12.2-3]{BroParPau19}.
We denote by $\Delta_x$ the unit Dirac mass at a point
$x$. 

\btheo\label{theo:quaterhyperbo} 
Let $\Ga,D^-,D^+$ be as above. There exists $\kappa>0$ such that, as
$s\ra+\infty$,
$$
\N_{D^-,\,D^+}(s)=c(D^-,D^+)\;
e^{(4n+2)\,s}\;\big(1+\operatorname{O}(e^{-\kappa s})\big)\;.
$$
Furthermore, the origins of the common perpendiculars from $D^-$ to the
images of $D^+$ under the elements of $\Ga$ equidistribute in
$\partial D^-$ to the induced Riemannian measure: as
$s\ra+\infty$,
\begin{equation}\label{eq:distribhorobcomplexhyp}
\frac{2\,(2n+1)\,\Vol(\Ga_{D^-}\bs D^-)}{c(D^-,D^+)}
\;e^{-(4n+2)\,s}\;\sum_{x\in\partial D^-} m_{s}(x) \;
\Delta_{x}\;\weakstar\; \vol_{\partial D^-}\,.
\;\;\;\Box
\end{equation}
\etheo

For smooth functions $\psi$ with compact support on $\partial D^-$,
there is an error term in the equidistribution claim of Theorem
\ref{theo:quaterhyperbo} when the measures on both sides are
evaluated on $\psi$, of the form $\bigO(e^{-\kappa s}\,\|\psi\|_\ell)$
where $\kappa>0$ and $\|\psi\|_\ell$ is the Sobolev norm of $\psi$ for
some $\ell\in\NN$.

\bigskip
We now apply Theorem \ref{theo:quaterhyperbo} in order to prove an
analog of Mertens's formula and Neville's equidistribution theorem in the
quaternionic Heisenberg group. See for example the Introduction of
\cite{ParPau17MA} for an explanation of the name.

Let $\mmm$ be a nonzero bilateral ideal in $\OOO$ stable by
conjugation.  As defined in Equations
\eqref{eq:heisintro}--\eqref{eq:heisaction} in the Introduction, the
action by shears on $\OOO\times \OOO \times \OOO$ of the nilpotent
group $\N(\OOO)$ preserves $\OOO\times\mmm\times\mmm$. We will study
the asymptotic of the counting function $\Psi_{\mmm}$, where, for
every $s\geq 0$, the number $\Psi_{\mmm}(s)$ is the cardinality of
$$
\N(\OOO)\bs
\big\{(a,\alpha,c)\in\OOO\times\mmm\times\mmm\;:\;
\tr(a\,\overline{c})=\n(\alpha),\;_\OOO\langle a,\alpha,c\rangle=\OOO,
\;0<\n(c)\leq s\big\}\,.
$$

We endow the ring $\OOO/\mmm$ with the involution induced by the
quaternionic conjugation. Let $\operatorname{U}_q(\OOO/\mmm)$ be the finite
group of $3\times3$ matrices in $\OOO/\mmm$, preserving the Hermitian
form $-\overline{z_0}\,z_2- \overline{z_2}\,z_0 + \overline{z_1}\,z_1$ on
$(\OOO/\mmm)^3$. Let $B_q(\OOO/\mmm)$ be its upper triangular
subgroup.

\btheo\label{theo:countHeis} 
There exists $\kappa>0$ such that, as $s\ra+\infty$,
$$
\Psi_{\mmm}(s)= 
\frac{204\,120\;D_A^{\;4}\;|B_q(\OOO/\mmm)|}
{\pi^8\;m_A\;|\OOO^\times|\;\prod_{p|D_A}(p-1)(p^2+1)(p^3-1)\;
|\operatorname{U}_q(\OOO/\mmm)|} 
\;s^5\;(1+\bigO(s^{-\kappa}))\,.
$$
\etheo

The particular case $\mmm=\OOO$ gives Theorem \ref{theo:countintro} in
the introduction. We will prove this result simultaneously with the
next one. We endow the Lie group $\HHeis_{7}$ with its Haar measure
$\haarhheis$ defined in the introduction. The following result is an
equidistribution result of the set of $\QQ$-points (satisfying some
congruence properties) in $\HHeis_{7}$, seen as the set of $\RR$-points
of a $\ZZ$-form of a $\QQ$-algebraic group with set of $\QQ$-points
$\HHeis_{7}\cap (A\times A)$ and set of $\ZZ$-points $\N(\OOO)$.  The
particular case $\mmm=\OOO$ gives Theorem \ref{theo:equidisintro} in
the introduction.

\btheo\label{theo:equidisHeis} As $s\ra+\infty$, we have
\begin{align*}
\frac{\pi^8\;m_A\;|\OOO^\times|\;\prod_{p|D_A}(p-1)(p^2+1)(p^3-1)\;
|\operatorname{U}_q(\OOO/\mmm)|}
{816\,480\;D_A^2\;|B_q(\OOO/\mmm)|}\;s^{-5}&\;\times \\
\;\sum_{\substack{(a,\,\alpha,\,c)\in \OOO\times\mmm\times\mmm,
\;0<\n(c)\leq s \\ \tr(a\,\overline{c})=\n(\alpha),
\;_\OOO\langle a,\,\alpha,\,c\rangle=\OOO}}\;
&\Delta_{(ac^{-1},\,\alpha c^{-1})}\;\weakstar\;\haarhheis\;.
\end{align*}
\etheo

As in Theorem \ref{theo:quaterhyperbo}, for smooth functions
$\psi$ with compact support on $\HHeis_7$, there is an error term in
this equidistribution result when the measures on both sides are
evaluated on $\psi$, of the form $\bigO(s^{-\kappa}\,\|\psi\|_\ell)$
where $\kappa>0$ and $\|\psi\|_\ell$ is the Sobolev norm of $\psi$ for
some $\ell\in\NN$.

\medskip
\noindent{\bf Proofs of Theorem \ref{theo:countHeis} and Theorem
  \ref{theo:equidisHeis}.}  We start by introducing the notation used
in these proofs.

We consider the quaternionic Hermitian form $q$ defined in Section
\ref{sec:quathypspace} with $n=2$.  For every subgroup $G$ of
$\operatorname{U}_q$, we denote by $\overline{G}$ its image in
$\PU_q$, and again by $g$ the image in $\PU_q$ of any element $g$ of
$\operatorname{U}_q$.

We consider the lattice $\Ga=\operatorname{U}_q(\OOO)$ in
$\operatorname{U}_q$ defined in Section \ref{sec:cuspnumb}, so that
$\overline{\Ga}=\PU_q(\OOO)$. We denote by $\Ga_\mmm$ the Hecke
congruence subgroup of $\Ga$ modulo $\mmm$, that is the preimage, by
the group morphism $\Ga\ra \operatorname{U}_q(\OOO/\mmm)$ of reduction
modulo $\mmm$, of the upper triangular subgroup $B_q(\OOO/\mmm)$.
Since $-\id\in\Ga_\mmm$, we have
\begin{equation}\label{eq:calcindex}
[\;\overline{\Ga}\,:\;\overline{\Ga_\mmm}\;]=[\Ga:\Ga_\mmm]=
\frac{|\operatorname{U}_q(\OOO/\mmm)|}{|B_q(\OOO/\mmm)|}\;.
\end{equation}

We denote by $\Ga_{\H_1}$ the stabiliser in $\Ga_\mmm$ of the horoball
$\H_1$ defined in Equation \eqref{eq:defhoroinfty}. It is equal to
$B_q\cap\Ga_\mmm$ where $B_q$ has been defined in Section
\ref{sec:quathypspace}, since an element of $\Gamma$ fixes $\infty$ if
and only if it preserves $\H_1$.  The group $\Ga_{\H_1}$ is
independent of $\mmm$, by the definition of $\Ga_\mmm$.

The projection map from $\Ga_{\H_1}$ to $\overline{\Ga_{\H_1}}$ is
$2$-to-$1$ since $-\id\in \Ga_{\H_1}$. We identify the lattice
$\N(\OOO)$ of $\HHeis_7$ with its image $\overline{\N(\OOO)}$ by the
embedding of $\HHeis_7$ in $\PU_q$ defined in Section
\ref{sec:quathypgeom}. The system of equations \eqref{eq:equationsUq}
gives
\begin{equation}\label{eq:heisindexinparab}
[\;\overline{\Ga_{\H_1}}:\overline{\N(\OOO)}\;]=
\frac{1}{2}\;[\Ga_{\H_1}:\N(\OOO)]=
\frac{1}{2}\;|\OOO^\times|^2\,.
\end{equation}
The following result gives in particular the computation of the volume
of the cusp at infinity for $\overline{\Ga}$.

\blemm \label{lem:equalhaar}
The Haar measure $\lambda_7$ on $\HHeis_{4n-1}$ defined in
Equation \eqref{eq:defhaarlambda} coincides with the Haar measure
$\haarhheis$ defined in the introduction, that is, the total mass of
the measure induced by $\lambda_7$ on $\N(\OOO)\bs\HHeis_{7}$ is
$\frac{D_A^2}{4}$. Furthermore
\begin{equation}\label{eq:volcusp}
\Vol(\;\overline{\Ga_{\H_1}}\;\bs\,\H_1)=
\frac{D_A^2}{160\,|\,\OOO^\times|^2}\;.
\end{equation}
\elemm

For instance, if $D_A=2$ and $\OOO$ is the Hurwitz order
$\frac{1+i+j+k}{2}\ZZ+\ZZ i+\ZZ j+\ZZ k$, which has $24$ units and
$\Im\,\OOO=\ZZ i+\ZZ j+\ZZ k$, then
$\vol(\overline{\Ga_{\H_1}}\bs\H_1) =\frac1{23040}$, to be compared
with \cite[Prop.5.8]{KimPar03} for a related computation.

\medskip
\dem Note that $\tr :\HH\ra\RR$ is a fibration, with fiber
$\frac{t}{2} +\Im\, \HH$ over $t \in \RR$. The Lebesgue measure of the
Euclidean space $\HH$ disintegrates by this fibration over the
Lebesgue measure of $\RR$, with conditional measures on the fiber
$1/2$ the Lebesgue measure of the fiber: $dx_0dx_1dx_2dx_3 =
(\frac{1}{2}dx_1dx_2dx_3) d(2x_0)$. Since the map $\tr$ is additive,
since it maps $\OOO$ onto $\ZZ$ with kernel $\Im \,\OOO$ as recalled
in Section \ref{sect:quaternion}, this implies that
$\Vol(\Im\,\OOO\,\bs \,\Im\,\HH) = 2\Vol(\OOO\,\bs\,\HH)$. Again by
the surjectivity of $\tr:\OOO\ra\ZZ$, for every $w\in\OOO$, the set
$\{w_0\in\OOO\;:\;\tr w_0 =\n(w)\}$ is a translate of
$\Im\,\OOO$. Hence by Equations \eqref{eq:horosphecoord} and
\eqref{eq:covolmaxorder}, we have
$$
\lambda_7(\N(\OOO)\,\bs\,\HHeis_7)= 
2\Vol(\Im\,\OOO\,\bs \,\Im\,\HH)\Vol(\OOO\,\bs\,\HH)=
4\Vol(\OOO\,\bs\,\HH)^2=\frac{D_A^{\;2}}{4}\;.
$$
This proves the first claim of Lemma \ref{lem:equalhaar}.

Now, by Equation \eqref{eq:heisindexinparab}, by Lemma
\ref{lem:horobvol}, by Equation \eqref{eq:volHun}, we have
\begin{align*}
\Vol(\;\overline{\Ga_{\H_1}}\;\bs\,\H_1)&=
\frac2{|\OOO^\times|^2}\Vol(\;\overline{\N(\OOO)}\;\bs\,\H_1)=
\frac1{5\;|\OOO^\times|^2}
\Vol(\;\overline{\N(\OOO)}\;\bs\,\partial \H_1)\\ & =
\frac1{40\;|\OOO^\times|^2}\;
\lambda_7(\N(\OOO)\,\bs\,\HHeis_7)=
\frac{D_A^{\;2}}{160\;|\OOO^\times|^2}\;. \;\;\;\Box
\end{align*}

\medskip
We need one more notation before giving the proof of Theorem
\ref{theo:countHeis}. Consider an element $g\in \Ga_\mmm$ such that
$g\H_1$ and $\H_1$ are disjoint (there are only finitely many double
classes $[g]\in \overline{\Ga_{\H_1}} \,\bs \overline{\Ga_\mmm}/
\,\overline{\Ga_{\H_1}}$ for which this is not the case). We denote by
$\ell(\delta_g)$ the length of the common perpendicular $ \delta_g$
between $g\,\H_1$ and $\H_1$. If $\begin{pmatrix} a_g \\ \alpha_g
  \\ c_g \end{pmatrix}$ is the first column of $g$, then $g\cdot
\infty = [a_g:\alpha_g:c_g]$.

We use the following facts in the system of equations below.

$\bullet$~ For the first equality, note that the cardinality of each
nonempty fiber of the projection map from $\big\{(a,\alpha,c)\in\OOO
\times \mmm\times\mmm\;:\;_\OOO\langle a,\alpha,c\rangle=\OOO\big\}$
to $\PP^2_{\rm r}(\HH)$ is $|\OOO^\times|$ and that the projection
from $\N(\OOO)$ to $\overline{\N(\OOO)}$ is injective.

$\bullet$~ The second and third equalities follow from Proposition
\ref{prop:identiforbit} (2).

$\bullet$~ The fourth equality follows from Lemma
\ref{lem:disthoroghoro}.

$\bullet$~ The fifth equality follows by Equation
\eqref{eq:heisindexinparab} and by the definition of the counting
function $\N_{\H_1,\,\H_1}$.

$\bullet$~ The sixth equality follows from the first claim in Theorem
\ref{theo:quaterhyperbo} with $n=2$, $\Ga=\overline{\Ga_\mmm}$ and
$D^-=D^+=\H_1$.

 $\bullet$~ The last equality follows from Equations
\eqref{eq:volcusp} and \eqref{eq:calcindex} and from Theorem
\ref{intro:computvol}.

\noindent  We hence have, for some $\kappa>0$ and for every $s>0$,
\begin{align*}
\Psi_\mmm(s)& = |\OOO^\times|\;\;\card\;\;_{\overline{\N(\OOO)}}\,\bs
\Big\{[a:\alpha:c]\in\PP^2_{\rm r}(\HH)\;:\;
\begin{array}{c}(a,\alpha,c)\in\OOO\times\mmm\times\mmm,\\
_\OOO\langle a,\alpha,c\rangle=\OOO,\\
\tr(a\,\overline{c})=\n(\alpha),\;0<\n(c)\leq s\end{array}
\Big\}\nonumber\\ 
&=
|\OOO^\times|\;\;\card\;\;_{\overline{\N(\OOO)}}\,\bs
\big\{[a:\alpha:c]\in\overline{\Ga_\mmm}\cdot\infty\;:\;
\begin{array}{c}(a,\alpha,c)\in\OOO\times\mmm\times\mmm,\\
 _\OOO \langle a,\alpha,c\rangle=\OOO,\;0<\n(c)\leq s\end{array}
\big\}\nonumber\\ 
&=
|\OOO^\times|\;\;\card\;
\big\{[g]\in\overline{\N(\OOO)}\,\bs\,
\overline{\Ga_\mmm}\,/\,\overline{\Ga_{\H_1}}\;:\;0<\n(c_g)\leq s
\big\}\nonumber\\ 
&=
|\OOO^\times|\;
[\,\overline{\Ga_{\H_1}}:\overline{\N(\OOO)}\,]\;\card
\big\{[g]\in\overline{\Ga_{\H_1}}\;\bs\,
\overline{\Ga_\mmm}\,/\,\overline{\Ga_{\H_1}}\;:\; \ell(\delta_g)\leq
\frac{\ln s}{2}-\ln 2 \big\}+\bigO(1)\nonumber\\ 
&=
\frac12\;|\OOO^\times|^3\;
\;\N_{\H_1,\,\H_1}\big(\,\frac{\ln s}{2}-\ln 2\big)
+\bigO(1) \nonumber \\ 
&=
\frac{15\;|\OOO^\times|^3\;
\big(\Vol(\;\overline{\Ga_{\H_1}}\;\bs\,\H_1)\big)^2}
{\pi^4\,\,\Vol(\;\overline{\Ga_\mmm}\;\bs\, \HH^2_\HH)} 
\;s^5\,(1+\bigO(s^{-\kappa}))\\
&=
\frac{204\,120\;D_A^{\;4}\;|B_q(\OOO/\mmm)|}
{\pi^8\;m_A\;|\OOO^\times|\;\prod_{p|D_A}(p-1)(p^2+1)(p^3-1)
\;|\operatorname{U}_q(\OOO/\mmm)|}
 \;s^5\,(1+\bigO(s^{-\kappa}))\;. 
\end{align*} 
This concludes the proof of Theorem \ref{theo:countHeis}.

\medskip
Let us prove now Theorem \ref{theo:equidisHeis}. The orthogonal
projection map $f:\partial_\infty\HH^2_\HH-\{\infty\}\ra \partial
\H_1$ is the homeomorphism defined by $[w_0:w:1]\mapsto (\zeta=w,
u=2\,\Im\,w_0,1)$ using the homogeneous coordinates on
$\partial_\infty\HH^2_\HH-\{\infty\}$ and the horospherical
coordinates on $\partial \H_1$ (see Equation
\eqref{eq:horosphecoord}). Let $x\in\partial \H_1$ be the origin of a
common perpendicular of length at most $s$ from $\H_1$ to a horoball
$\ga\H_1$ for some $\ga\in\Ga_\mmm$ not fixing $\infty$. The point $x$
is the orthogonal projection on $\H_1$ of the point at infinity of
this horoball $\ga\H_1$. This point at infinity may be written
$[ac^{-1},\alpha c^{-1} :1]$ for some triple $(a,\,\alpha,\,c) \in
\OOO\times \mmm \times\mmm$ with $_\OOO\langle a,\,\alpha,\,c\rangle
=\OOO$, $\tr(a\,\overline{c})=\n(\alpha)$ and $0<\n(c)\leq 4\,e^{2s}$
(using Lemma \ref{lem:disthoroghoro}). Such a writing is not unique,
there are exactly $|\OOO^\times|$ such triples. Hence by the second
claim of Theorem \ref{theo:quaterhyperbo} with $\Ga=
\overline{\Ga_\mmm}$ and $D^-=D^+=\H_1$, using the horospherical
coordinates on $\partial\H_1$, we have, as $s\ra+\infty$,
\begin{equation}\label{eq:appliPP17ETDS}
\frac{\pi^4\;\Vol(\,\overline{\Ga_\mmm}\,\bs \HH^2_\HH)}
{3072\;\Vol(\,\overline{\Ga_{\H_1}}\,\bs\H_1)\;|\OOO^\times|}\;e^{-10\,s}
\;\sum_{\substack{(a,\,\alpha,\,c)\in \OOO\times\mmm\times\mmm,
\;0<\n(c)\leq 4\,e^{2s} \\ \tr(a\,\overline{c})=\n(\alpha),
\;_\OOO\langle a,\,\alpha,\,c\rangle=\OOO}}\;
\Delta_{(ac^{-1},\,2\,\Im(\alpha c^{-1}),\,1)}\;\weakstar\; \vol_{\partial \H_1}\;.
\end{equation}

Recall that the Haar measure $\lambda_7$ on
$\HHeis_7=\partial_\infty\HH^2_\HH-\{\infty\}$, defined in Equation
\eqref{eq:defhaarlambda}, coincides with the Haar measure $\haarhheis$
by Lemma \ref{lem:equalhaar}. Its image by the above map $f$ is, by
Equation \eqref{eq:volhorosphzetau},
$$
f_*\haarheis=f_*\,\lambda_7=8\vol_{\partial \H_1}\;.
$$
Using the change of variables $s\mapsto 4\,e^{2s}$ and the continuity of
the pushforward by $f^{-1}$ of the measures on $\partial \H_1$ applied
to Equation \eqref{eq:appliPP17ETDS}, we hence have, as $s\ra
+\infty$,
$$
\frac{8\;\pi^4\;\Vol(\,\overline{\Ga_\mmm}\,\bs \HH^2_\HH)}
{3\;|\OOO^\times|\;\Vol(\,\overline{\Ga_{\H_1}}\,\bs\H_1)}\;s^{-5}
\;\sum_{\substack{(a,\,\alpha,\,c)\in \OOO\times\mmm\times\mmm,
\;0<\n(c)\leq s \\ \tr(a\,\overline{c})=\n(\alpha),
\;_\OOO\langle a,\,\alpha,\,c\rangle=\OOO}}\;
\Delta_{(ac^{-1},\,\alpha c^{-1})}\;\weakstar\;\haarheis\;.
$$
Finally, Theorem \ref{theo:equidisHeis} follows from this, from Equations
\eqref{eq:volcusp} and \eqref{eq:calcindex} and from Theorem
\ref{intro:computvol}.
\cqfd

\medskip 
\brema\label{rem:highdim} {\rm Theorems \ref{theo:countHeis} and
 \ref{theo:equidisHeis} have generalisations in higher
dimension.  Theorem~\ref{theo:quaterhyperbo} (which is valid in any
dimension), applied with $\Ga=\PU_q(\OOO)$ and with $D^-=D^+$ the
horoball of points in $\hnh$ with last horospherical coordinates at
least $1$, gives a counting and equidistribution result  of the orbit $\Ga\cdot \infty-
\{\infty\}$ in $\partial_\infty\HHeis_{4n-1}$with error
term. The volume of $\Ga\bs \hnh$ could be computed using
\cite{EmeKim18}, up to computing the index of $\Ga$ in a principal
arithmetic subgroup containing it. The volume of the cusp
corresponding to $\infty$ in $\Ga\bs \hnh$ may also be computed by the
same method as for the proof of Equation \eqref{eq:volcusp}.

Other counting and equidistribution results of arithmetically defined
points in the quaternionic Heisenberg group $\HHeis_{4n-1}$ may be
obtained by varying the cusp (when $n=2$ and $h_A\neq 1$, there are at
least two cusps by Theorem \ref{intro:cuspnumb}), the integral
quaternionic Hermitian form $q$ of Witt signature $(1,n)$ and the
arithmetic lattice $\Ga$ in $\operatorname{U}_q$. }  \erema

{\small \bibliography{/users/jouniparkkonen/Documents/Latex/viitteet.bib} }

\bigskip
{\small
\noindent \begin{tabular}{l} 
Department of Mathematics and Statistics, P.O. Box 35\\ 
40014 University of Jyv\"askyl\"a, FINLAND.\\
{\it e-mail: jouni.t.parkkonen@jyu.fi}
\end{tabular}
\medskip

\noindent \begin{tabular}{l}
Laboratoire de math\'ematique d'Orsay, UMR 8628 Univ.~Paris-Sud, CNRS\\
Universit\'e Paris-Saclay,
91405 ORSAY Cedex, FRANCE\\
{\it e-mail: frederic.paulin@math.u-psud.fr}
\end{tabular}
}

\end{document}